\documentclass[12pt, reqno, twoside]{amsart}
\usepackage{amsmath, amsthm, mathrsfs, amscd, amsfonts, amssymb, graphicx, color}
\usepackage[bookmarksnumbered, colorlinks, plainpages]{hyperref}
\hypersetup{colorlinks=true,linkcolor=red, anchorcolor=green, citecolor=cyan, urlcolor=red, filecolor=magenta, pdftoolbar=true}
\usepackage{enumitem}
\usepackage{enumitem}
\newtheorem{theorem}{Theorem}[section]
\newtheorem{lemma}[theorem]{Lemma}

\newtheorem{corollary}[theorem]{Corollary}
\theoremstyle{definition}

\theoremstyle{remark}
\newtheorem{remark}[theorem]{Remark}
\numberwithin{equation}{section}

\begin{document}

\setcounter{page}{1}

\title[Gradient estimates for the porous medium equation]
{The nonlinear porous medium equation for the $f$-Laplacian: Hamilton-Souplet-Zhang type gradient estimates and implications}

\author[A. Taheri]{Ali Taheri}

\author[V. Vahidifar]{Vahideh Vahidifar}

\address{School of Mathematical and  Physical Sciences, 
University of Sussex, Falmer, Brighton, United Kingdom.}
\email{\textcolor[rgb]{0.00,0.00,0.84}{a.taheri@sussex.ac.uk}} 
 
\address{School of Mathematical and  Physical Sciences, 
University of Sussex, Falmer, Brighton, United Kingdom.}
\email{\textcolor[rgb]{0.00,0.00,0.84}{v.vahidifar@sussex.ac.uk}}

\subjclass[2010]{53C44, 58J60, 58J35, 60J60}

\keywords{Smooth metric measure spaces, $f$-Laplacian, Nonlinear diffusion, Porous medium equation, 
Super Perelman-Ricci flow, Hamilton-Souplet-Zhang estimates, Ancient solutions}

%\date{November 4, 2019}

\begin{abstract}
This article presents new gradient estimates for positive solutions to the nonlinear porous medium equation (NPME) 
in the context of smooth metric measure spaces. The diffusion operator here is the $f$-Laplacian and the gradient 
estimates of interest are mainly of Hamilton-Souplet-Zhang types. These estimates are established using a variety of methods 
and techniques and several implications, most notably, to parabolic Liouville-type results and characterisation of ancient 
solutions are given. The problem is posed in the general framework where the metric and potential evolve with 
time and the proofs make use of natural lower bounds on the time derivative of the metric and the Bakry-\'Emery 
$m$-Ricci curvature tensors. Our results extend and improve various existing ones in the literature. 
\end{abstract}

\maketitle
 
{ 
\hypersetup{linkcolor=black}
\tableofcontents
}

\allowdisplaybreaks

\section{Introduction} \label{sec1}

In this paper we are primarily concerned with nonlinear diffusion equations of porous medium type, or equivalently, 
slow diffusion type, in the framework of smooth metric measure spaces (hereafter we write SMMS for brevity). The motivation 
for this study is prompted by the recent surge of interest in understanding the analytic and dynamic behaviour of 
solutions, the probabilistic features of diffusion on SMMS and the role that geometry and curvature on the one 
hand and nonlinearity on the other play in forming such behaviour and features 
({\it see}, e.g. \cite{AM, Bak, BBGM, BDM, Bon, DaskKe, Otto, Pe02, VC, Wang, Zhang}).

Towards this end, let $(\mathscr M,g)$ be a complete Riemannian manifold of dimension $n \ge 2$, $d\mu=e^{-f} dv_g$ be 
a positive weighted measure associated with the metric $g$ and potential $f$ and $dv_g$ be the standard Riemannain 
volume measure on $\mathscr M$. The triple $(\mathscr M, g, d\mu)$ is referred to as a smooth metric measure space (or a weighted manifold). 
Our prime interest in this paper is the nonlinear porous medium equation associated with the triple $(\mathscr M,g,d\mu)$ 
that can be written for a positive $u=u(x,t)$ in the form:
\begin{align} \label{eq11}
\square_p u(x,t) = \partial_t u (x,t) - \Delta_f u^p (x,t) = \mathscr N(t,x,u(x,t)), \qquad p>1.  
\end{align}
The operator $\Delta_f$ in \eqref{eq11} is the $f$-Laplacian (also called Witten or weighted Laplacian) 
that acts on functions $w \in \mathscr{C}^2(\mathscr M)$ as  
\begin{equation} \label{f-Lap-definition}
\Delta_f w = e^f {\rm div}(e^{-f} \nabla w) = \Delta w - \langle \nabla f, \nabla w\rangle.
\end{equation} 
The exponent $p>1$ is a fixed constant 
and ${\mathscr N}={\mathscr N}(t,x,u)$ is a sufficiently smooth {\it nonlinearity} or forcing term that depends 
on the space-time variables $(x, t)$ and the dependent variable $u>0$. The structural 
features of this nonlinearity have analytical and physical significances and one of our aims is 
to understand the ways it influences the estimates and subsequently the solutions. Our task 
is to develop gradient estimates of elliptic and Hamilton-Souplet-Zhang 
types for the positive solutions $u=u(x,t)$ to \eqref{eq11} where we use a variety of techniques 
and ideas to deal with different exponent ranges and investigate some of their important implications.

The porous medium equation (PME) arises in many contexts with many applications. Indeed for various values of $p>1$ 
it arises in different models of diffusive phenomenon: Boussinesq's model of ground water infiltration, flow of gas 
in porous media, heat radiation in plasma, liquid thin films and so on. We refer the interested reader to the 
monographs and texts \cite{Aron, PLu, DaskKe, Vaz} and the references therein for a comprehensive coverage of these 
applications and the underlying theory. The PME can be seen as the simplest generalisation 
of the linear heat equation to the nonlinear context (coinciding with the latter when $p=1$). The use of terminology 
{\it slow diffusion} here is intertwined with the exponent range $p>1$ in \eqref{eq11}, where, due to the density dependent 
diffusion coefficient vanishing at $u=0$ -- as can be seen by writing $\Delta_f u^p=pe^f {\rm div} (e^{-f} u^{p-1} \nabla u)$ 
-- disturbances of $u=0$ propagate forward in time with finite speed. This is in sharp contrast to the linear heat as well as 
the fast (and very fast) diffusion case $0<p<1$ that has a completely different nature and morphology and will be treated in a separate paper.

Following standard theory ({\it see} \cite{Aron, DaskKe, Vaz}), we consider in the SMMS context under study here, the second 
order linear (space-time dependent variable coefficient) evolution operator defined by 
\begin{align} \label{Lpv-definition-introduction}
\mathscr L_v^p = \partial_t -(p-1) v(x,t) \Delta_f, 
\end{align}
where $v$ relates to $u$ via the transformation $v =pu^{p-1}/(p-1)$ (the pressure transform) 
and $p>1$. This operator plays an important role throughout the paper, particularly, in that, if $u$ is a 
positive solution to \eqref{eq11}, then the pressure $v$ is a positive solution to a related equation involving the 
linear (heat-type) operator $\mathscr L_v^p$. See Sections \ref{sec3} and \ref{sec7} for more and the 
discussion surrounding \eqref{SPR-p-substitute-intro} at the end of this section.

The $f$-Laplacian $\Delta_f$ is a symmetric diffusion operator with respect to the invariant weighted measure 
$d\mu$ and arises in a variety of contexts ranging from probability theory and geometry to quantum field theory 
and statistical mechanics \cite{Bak, Wang}. 
It is a natural generalisation of the Laplace-Beltrami 
operator to the SMMS context and it coincides with the latter precisely when the potential $f$ is spatially constant. 
By an application of the integration by parts formula it can be seen that for 
$u, v \in {\mathscr C}_0^\infty(\mathscr M)$ it holds  
\begin{equation}
\int_\mathscr M u \Delta_f v \, d\mu 
= - \int_\mathscr M \langle \nabla u, \nabla v \rangle \, d\mu 
= \int_\mathscr M v \Delta_f u \, d\mu.
\end{equation}

As for the geometry and curvature properties of the triple $(\mathscr M,g,d\mu)$ 
we have a hierarchy of second order symmetric tensor fields on $\mathscr M$ defined by,  
\begin{equation} \label{Ricci-m-f-intro}
{\mathscr Ric}^m_f(g) := {\mathscr Ric}(g) + \nabla\nabla f - \frac{\nabla f \otimes \nabla f}{m-n}, \qquad m \ge n, 
\end{equation}
called the Bakry-\`Emery $m$-Ricci curvature tensors. Here ${\mathscr Ric}(g)$ is the usual Riemannain Ricci 
curvature tensor of $g$, $\nabla \nabla f={\rm Hess}(f)$ denotes the Hessian of $f$, and $m \ge n$ is a constant 
({\it see} \cite{Bak}). For the sake of clarification, in relation to \eqref{Ricci-m-f-intro}, when 
$m=n$, by convention $f$ is only allowed to be a constant, thus giving ${\mathscr Ric}^m_f(g)={\mathscr Ric}(g)$. 
We also allow for $m = \infty$ in which case by formally passing to the limit in \eqref{Ricci-m-f-intro} we define, 
\begin{equation} \label{Ricf def eq} 
{\mathscr Ric}_f(g) = {\mathscr Ric}(g) + \nabla\nabla f := {\mathscr Ric}^\infty_f(g).
\end{equation}

By an easy inspection it is seen that having a lower bound on ${\mathscr Ric}_f^m(g)$ is a stronger condition than having 
the same lower bound on ${\mathscr Ric}_f(g)$, that is,  
\begin{equation}
{\mathscr Ric}_f^m(g) \ge {\mathsf k} g \implies {\mathscr Ric}_f(g) \ge {\mathsf k} g,
\end{equation} 
but the reverse implication is not true (here ${\mathsf k} \in {\mathbb R}$, $m \ge n$ are finite constants). 
This remark is useful when comparing the results for nonlinear diffusion equations here with 
the ones for the linear heat-type equations. %The latter uses lower bounds on ${\mathscr Ric}_f(g)$ whereas the 
%former required lower bounds on ${\mathscr Ric}_f^m(g)$. 

The crucial identity relating the $f$-Laplacian $\Delta_f$ to the Bakry-Emery Ricci curvature tensor ${\mathscr Ric}_f(g)$ 
in this context is the weighted Bochner-Weitzenb\"ock formula (generalising, in turn, the classical version of the identity to 
the context of SMMS, {\it see}, e.g., \cite{Aub, [LiP], Lott}) asserting that or any function $w$ of class ${\mathscr C}^3(\mathscr M)$ 
we have, 
\begin{equation} \label{Bochner-1}
\frac{1}{2} \Delta_f |\nabla w|^2 - \langle \nabla w, \nabla \Delta_f w \rangle
= |\nabla\nabla w|^2 + {\mathscr Ric}_f (\nabla w, \nabla w).
\end{equation}

While this identity is not immediately applicable to the Bakry-\'Emery $m$-Ricci curvature tensor, upon adding and 
subtracting the symmetric second order rank-one tensor $[\nabla f \otimes \nabla f]/(m-n)$ to \eqref{Bochner-1}, 
making use of \eqref{Ricci-m-f-intro} and then an application of the Cauchy-Schwarz inequality, one can arrive 
at a counterpart {\it inequality} in the form 
\begin{align} \label{Bochner-2}
\frac{1}{2} \Delta_f |\nabla w|^2 - \langle \nabla w, \nabla \Delta_f w \rangle
&\ge \frac{(\Delta w)^2}{n} + \frac{[\nabla f \otimes \nabla f]}{m-n} (\nabla w, \nabla w) 
+ {\mathscr Ric}^m_f (\nabla w, \nabla w) \nonumber \\
&\ge \frac{ (\Delta_f w)^2}{m} + {\mathscr Ric}^m_f (\nabla w, \nabla w).
\end{align}

Interestingly from \eqref{Bochner-1} or \eqref{Bochner-2} it follows that subject to a curvature lower bound 
in the form ${\mathscr Ric}_f(g) \ge {\mathsf k} g$ or ${\mathscr Ric}^m_f(g) \ge {\mathsf k} g$ with 
${\mathsf k} \in {\mathbb R}$, the diffusion operator $\Delta_f$ satisfies the curvature-dimension 
condition ${\rm CD}({\mathsf k}, \infty)$ or ${\rm CD}({\mathsf k}, m)$ 
respectively which in turn plays an important role in the study of the dynamics, convergence to equilibrium 
and the derivation various functional and geometric inequalities associated with $\Delta_f$ ({\it see}, e.g., 
\cite{Bak} for more).

Gradient estimates lie at the core of geometric analysis and have huge applications. 
For heat and Schr\"odinger type equations on static manifolds they were first established 
in the seminal paper of Li and Yau \cite{[LY86]}. In the nonlinear setting, equations of heat-type 
with a logarithmic nonlinearity were among the first to be considered ({\it see} \cite{LiJ91, Ma}). 
In the context of SMMS the nonlinear diffusion version of these equations can be written in the 
form \eqref{eq11} with ${\mathscr N}(t,x,u)={\mathsf A}(t,x) u \log u$, specifically, 
\begin{equation} \label{eq1.4}
\square_p u(x,t) = \frac{\partial u}{\partial t} (x,t) - \Delta_f u^p (x,t) = {\mathsf A}(t,x) u(x,t) \log u(x,t). 
\end{equation}

The interest in such class of problems come from the natural links with geometric and functional inequalities on 
manifolds ({\it cf.} \cite{Bak, Gross, VC, Wang, Zhang} and the references therein for more). 
Other classes of equations closely relating to \eqref{eq1.4} include 
\begin{equation}
\square_p u(x,t) = \frac{\partial u}{\partial t} (x,t) - \Delta_f u^p (x,t) 
= {\mathsf A}(t,x) u^\alpha(x,t) \Gamma(\log u(x,t)) + {\mathsf B}(t,x) u^\beta(x,t), 
\end{equation}
with ${\mathscr N}(t,x,u) = {\mathsf A}(t,x) u^\alpha \Gamma(\log u) + {\mathsf B}(t,x) u^\beta$ 
where $\alpha$, $\beta$ are real exponents, ${\mathsf A}$, ${\mathsf B}$ are sufficiently smooth space-time functions and 
$\Gamma \in {\mathscr C}^1(\mathbb{R}, \mathbb{R})$. 
Another class of equations that have been extensively studied are Yamabe type equations. In the SMMS context, 
the nonlinear diffusion version of these equations take the form 
\begin{align} \label{eq11s}
\square_p u(x,t) = \frac{\partial u}{\partial t} (x,t) - \Delta_f u^p (x,t) = {\mathsf A}(t,x) u^\alpha(x,t) + {\mathsf B}(t,x) u^\beta(x,t). 
\end{align}
A far reaching generalisation of \eqref{eq11s} with a superposition of power-like nonlinearities consist of equations in the form 
\begin{align}
\square_p u(x,t) &= \frac{\partial u}{\partial t} (x,t) - \Delta_f u^p(x,t) \nonumber \\
&=  \sum_{j=1}^d {\mathsf A}_j(t,x) u^{\alpha_j} (x,t) + \sum_{j=1}^d {\mathsf B}_j(t,x) u^{\beta_j} (x,t). 
\end{align}
Here by comparison with \eqref{eq11} it is seen that the nonlinearity takes the form 
\begin{align}
{\mathscr N}(t,x,u) = \sum_{j=1}^d {\mathsf A}_j(t,x) u^{\alpha_j} + \sum_{j=1}^d {\mathsf B}_j(t,x) u^{\beta_j}, 
\end{align}
where ${\mathsf A}_j$, ${\mathsf B}_j$ (with $1 \le j \le d$) are sufficiently smooth space-time dependent coefficients and 
$\alpha_j \ge 0$, $\beta_j \le 0$ are real exponents ({\it see} \cite{Taheri-GE-1, Taheri-GE-2, TVahNA}).

The context in which we consider equation \eqref{eq11} is one where the metric tensor $g$ and the potential $f$ are 
allowed to be time dependent. Such problems have particularly moved to the forefront of research in geometric analysis since the work 
of G.~Perelman \cite{Pe02} and the study of forward and backward heat-type equations and more generally gradient flows on 
manifolds evolving under (super) Ricci or (super) Perelman-Ricci flows 
({\it see} also \cite{BaCP, Chow, SZ, Taheri-GE-1, Taheri-GE-2, TVahNA, TVDiffHar, Zhang}).

In this time dependent setting, the measure $d\mu = e^{-f} dv_g$, the Bakry-\'Emery $m$-Ricci curvature tensor 
${\mathscr Ric}_f^m(g)$ and the usual (metric dependent) differential operators $\nabla$, ${\rm div}, \Delta$ and 
$\Delta_f$ are all time dependent too which makes the analysis more challenging. 
In the course of the paper, it will become apparent that the super flow inequality  
\begin{align} \label{SPR-p-substitute-intro}
\frac{1}{2} \dfrac{\partial g}{\partial t} (x,t) + p u^{p-1}(x,t) {\mathscr Ric}^m_f(g)(x,t) \ge - \mathsf{k}g(x,t), 
\end{align}
relating the metric $g$, the potential $f$, and the positive solution $u$ to \eqref{eq11} with $p>1$ 
will play an interesting role. Here we note that  
\begin{align}
{\mathscr Ric}_f^m(g)(x,t) = {\mathscr Ric} (g)(x,t) + \nabla_{g} \nabla_{g} f (x,t) 
- [\nabla_{g} f \otimes \nabla_{g} f](x,t)/(m-n), 
\end{align}
(and we use the convention described earlier for $m=n$). Interestingly, in this nonlinear slow diffusion 
context, the super flow inequality \eqref{SPR-p-substitute-intro} presents itself as 
the natural substitute for the super Perelman-Ricci flow inequality that links to linear heat-type equations where $p=1$ 
(see \cite{Taheri-GE-1, Taheri-GE-2, TVahNA, TVDiffHar}). As will become evident later, the form of this super flow inequality 
is closely connected with the evolution operator ${\mathscr L}_p^v$ defined in \eqref{Lpv-definition-introduction}. 
(We refer the reader to Lemma \ref{Lem-8.2}, Lemmas \ref{Lem-3.1-new} and \ref{Lem-3.1-new-after} and the entire 
Section \ref{sec10} for more discussion on this.)

For the sake of establishing the main results and (local) estimates in this paper we shall make use of the lower bounds 
\begin{equation} \label{h-k-bounds-intro-discussion}
\partial_t g (x,t) \ge - 2h g(x,t), \qquad  {\mathscr Ric}^m_f(g)(x,t) \ge - k g(x,t), 
\end{equation}
for suitable constants $h$, $k \ge 0$ in a fixed compact space-time cylinder $Q_{R,T}$ with upper base centered 
at the reference point where the estimate is sought. These bounds on the one hand allow for the use of weighted 
(Laplace) comparison theorems and on the other provide control on the time derivative of geodesic distances 
which are two important ingredients in localisation and the proof of gradient estimates. The constants $h$, 
$k \ge 0$ will appear along with other bounds in the geometry-dependent terms in the ultimate formulation 
of the local estimates. In the static case, where the metric and potential do not depend on time 
(i.e., $\partial_t g \equiv 0$ and $\partial_t f \equiv 0$) we can take $h=0$ in \eqref{h-k-bounds-intro-discussion} 
and our results are immediately seen to cover this special case with the usual assumption of Ricci curvature lower 
bound (space only) of the type ${\mathscr Ric}_f^m(g) \ge - {\mathsf k} g$. It is clear that 
\eqref{h-k-bounds-intro-discussion} with a local upper bound on the positive solution $u=u(x,t)$ 
will immediately lead to \eqref{SPR-p-substitute-intro} with any choice of constant 
${\mathsf k} \ge h + p [\sup u]^{p-1} k \ge 0$, in view of,  
\begin{align}
\partial_t g (x,t) + 2p u^{p-1}(x,t) {\mathscr Ric}^m_f(g)(x,t) &\ge -2h g(x,t) - 2k pu^{p-1}(x,t) g(x,t) \\
&\ge -2 \{ h +  kp [\sup u]^{p-1} \} g(x,t) \ge - 2{\mathsf k} g(x,t). \nonumber 
\end{align}

\qquad \\
{\bf Plan of the paper.} The main estimates are presented in Theorem \ref{thm1} in Section \ref{sec2} and 
Theorem \ref{thm1-EQ-PME} in Section \ref{sec6} where both estimates are presented in their local forms. 
These estimates cover different exponent ranges and their proofs rely on different circle of ideas. The global 
form of these estimates follow by imposing appropriate global bounds on the solution, metric, potential as 
well as the Bakry-\'Emery $m$-Ricci curvature tensor and is presented in Theorem \ref{thm1-global} and Theorem 
\ref{thm1-EQ-PME-global} respectively. The static case (time independent metrics $g$ and potentials $f$) which 
is an important case is treated as a by-product of the above estimates and is presented in Theorem \ref{thm1-static} 
and Theorem \ref{thm1-EQ-PME-static} respectively. As a nice and useful application of these we are able 
to present parabolic Liouville results which lead to a characterisation of ancient solutions to \eqref{eq11}. 
These appear in turn in Theorem \ref{ancient-1} and Theorem \ref{ancient-2}. As such Sections \ref{sec2} 
and \ref{sec6} contain the main results of the paper whilst the remaining sections are devoted to 
developing the necessary apparatus and tools for establishing the results and proofs.

\qquad \\
{\bf Notation.} We write $z=z_+ + z_-$ with $z_+=\max(z, 0)$ and $z_-=\min(z, 0)$. 
Fixing a reference point $x_0 \in \mathscr M$ we denote by $d=d(x,x_0, t)$ the Riemannian distance between $x$ 
and $x_0$ on $\mathscr M$ with respect to $g=g(t)$. We write $\varrho=\varrho(x,x_0,t)$ for the geodesic radial 
variable measuring distance between $x$ and origin $x_0$ at time $t>0$. For a fixed space-time 
reference point $(x_0,t_0)$ and $R>0$, $T>0$ we define the space-time cylinder 
\begin{equation}
Q_{R,T} (x_0,t_0) \equiv \{ (x, t) | d(x, x_0, t) \le R, t_0-T \le t \le t_0 \} \subset \mathscr M \times [t_0-T, t_0].
\end{equation}
When the metric $g$ is time independent, we denote by $\mathscr{B}_\varrho(x_0) \subset \mathscr M$ the 
geodesic ball of radius $\varrho>0$ centred at $x_0$. It is evident that in this case we have 
\begin{equation}
Q_{R,T} (x_0,t_0) = \mathscr{B}_R(x_0) \times [t_0-T, t_0] \subset \mathscr M \times [t_0-T, t_0].
\end{equation}
Note that when the choice of the reference point $(x_0, t_0)$ 
is clear from the context we often abbreviate and write $d(x, t)$, $\varrho(x,t)$ or ${\mathscr B}_\varrho$, 
$Q_{R,T}$ respectively.

For a function of multiple variables, we denote its partial derivatives with subscripts accordingly. For instance 
for  $\Gamma=\Gamma(x,u)$ we denote its partial derivatives with respect to $x=(x_1, \dots, x_n)$ or $u$ by  
$\Gamma_x$ or $\Gamma_u$ respectively. Moreover we reserve the notation $\Gamma^x$ for 
the function $x \mapsto \Gamma (x, u)$ obtained by freezing the argument $u$ and viewing it as 
a function of $x$. In the sequel we frequently make use of the notations $\nabla \Gamma^x$ and $\Delta_f \Gamma^x$.

The Riemann, Ricci and Bakry-\'Emery $m$-Ricci curvature tensors associated 
with the metric $g$ and potential $f$, in local coordinates $(x_1, \dots, x_n)$, are given by:  
\begin{equation} \label{Eq-1.8}
[{\rm Riem}(g)]^\ell_{ijk} = \frac{\partial \Gamma^\ell_{jk}}{\partial x_i} - \frac{\partial \Gamma^\ell_{ik}}{\partial x_j} 
+ \Gamma^p_{jk} \Gamma^\ell_{ip}  - \Gamma^p_{ik} \Gamma^\ell_{jp},   
\end{equation}
\begin{equation} \label{Eq-1.9}
[{\mathscr Ric}(g)]_{ij} = \frac{\partial \Gamma_{ij}^k}{\partial x_k} - \frac{\partial \Gamma_{\ell j}^\ell}{\partial x_i} 
+ \Gamma^k_{ij} \Gamma^\ell_{\ell k} - \Gamma^\ell_{ik} \Gamma^k_{\ell j},  
\end{equation}
\begin{equation} \label{Eq-1.9-fm}
[{\mathscr Ric}_f^m(g)]_{ij} = \frac{\partial \Gamma_{ij}^k}{\partial x_k} - \frac{\partial \Gamma_{\ell j}^\ell}{\partial x_i} 
+ \Gamma^k_{ij} \Gamma^\ell_{\ell k} - \Gamma^\ell_{ik} \Gamma^k_{\ell j} 
+ \frac{\partial^2 f}{\partial x_i \partial x_j} - \Gamma^k_{ij} \frac{\partial f}{\partial x_k}
- \frac{1}{m-n} \frac{\partial f}{\partial x_i} \frac{\partial f}{\partial x_j}. 
\end{equation}
Likewise the $f$-Laplacian \eqref{f-Lap-definition} associated with the metric $g$ and potential $f$ takes the form: 
\begin{equation} \label{f-Lap-local}
\Delta_f = \underbrace{\frac{1}{\sqrt{|g|}} \frac{\partial}{\partial x_i} \left( \sqrt{|g|} g^{ij} \frac{\partial}{\partial x_j} \right)}_{\Delta} 
- \underbrace{g^{ij} \frac{\partial f}{\partial x_i} \frac{\partial}{\partial x_j}}_{\langle \nabla f, \cdot \rangle}.
\end{equation}
In \eqref{Eq-1.8}-\eqref{Eq-1.9-fm} above 
$\Gamma^k_{ij} = (g^{k\ell}/2) (\partial g_{j \ell}/\partial x_i + \partial g_{i \ell}/\partial x_j - \partial g_{ij}/\partial x_\ell)$, 
are the Christoffel symbols of $g$ whilst $g_{ij}$, $|g|$, $g^{ij} = (g^{-1})_{ij}$ are the components, 
determinant and the components of the inverse of $g$.

\section {A Hamilton-Souplet-Zhang estimate for \eqref{eq11}: $1< p < 1+1/(\sqrt{2m} +1)$}
\label{sec2}

In this section we present the first set of gradient estimates for the positive solutions to equation \eqref{eq11}. 
The estimates as will be seen below are valid in the exponent range 
$1<p<1+1/(\sqrt{2m} +1)$. Note that $m \ge n$ inherently has the role of dimension for the SMMS. 
Here the curvature lower bounds are expressed as ${\mathscr Ric}_f^m(g) \ge -(m-1)k g$ where 
the presence of $m$ is explicitly felt. Note that $m$ is not necessarily an integer.
\begin{theorem} \label{thm1}
Let $(\mathscr M, g, d\mu)$ be a complete smooth metric measure space with $d\mu=e^{-f} dv_g$. 
Assume the metric and potential are time dependent, of class $\mathscr{C}^2$ and 
that for suitable constants $k, h \ge 0$ and $m \ge n$ satisfy 
${\mathscr Ric}_f^m (g) \ge -(m-1)k g$, $\partial_t g \ge -2 h g$ 
in the space-time cylinder $Q_{R,T}$ with $R, T >0$. Let $u$ be a positive solution 
to \eqref{eq11} with $1 < p <1+1/(\sqrt{2m} +1)$ and $v =pu^p/(p-1)$ and 
$M=\sup_{Q_{R,T}} v$. Then there exists $C=C(p,\beta,m)>0$ such 
that for every $(x,t)$ in $Q_{R/2,T}$ with $t>t_0-T$ we have 
\begin{align}\label{eq-2.1}
\frac{|\nabla v|}{v^{\beta/2}}(x,t) 
\le C \left \{ \begin {array}{ll}
\sqrt h M^{(1-\beta)/2} + \left[ \dfrac{k^{1/4}}{\sqrt R} + \dfrac{1}{R} + \sqrt k \right] M^{1-\beta/2}
\\
\\
+ \dfrac{M^{(1-\beta)/2}}{\sqrt {t-t_0+T}} 
+ \sup_{Q_{R, T}} \left\{\left[\dfrac{|\Sigma_x(t,x,v)|}{v^{(3\beta-2)/2}}\right]^{1/3}\right\} 
\\
\\
+ \sup_{Q_{R, T}}\left\{v^{(1-\beta)/2} \left[2 \Sigma_v(t,x,v)-\dfrac{ \beta \Sigma(t,x,v)}{v} \right]_+^{1/2}\right\}
\end{array}
\right\}.
\end{align}
Here $\beta \in (\beta_1, \beta_2)$ where $\beta_1<\beta_2< 0$
are the roots of $\beta^2 + (2-p)/(p -1)\beta +m/2 =0$ and $\Sigma(t,x,v)$ is defined by
\begin{align}\label{EQ-eq-2.4}
\Sigma(t,x,v)= p [(p-1)v/p ]^{\frac{p -2}{p-1}} 
\mathscr N \left(t,x,[(p-1)v/p]^{\frac{1}{p-1}}\right).
\end {align}
\end{theorem}

\begin{remark}
The particular choice $\beta = -(2-p)/[2(p-1)]$ (easily seen to lie between the roots $\beta_1$, $\beta_2$) leads to 
the following formulation of the estimate \eqref{eq-2.1} that is recorded for future reference:
\begin{align}\label{eq-2.2}
v^{\frac{1}{4}\frac{2-p}{p-1}}|\nabla v|
\le C \left \{ \begin {array}{ll}
\sqrt h M^{\frac{1}{2}+\frac{1}{4}\frac{2-p}{p-1}} + \left[ \dfrac{k^{1/4}}{\sqrt R} + \dfrac{1}{R} + \sqrt k \right] M^{1+\frac{1}{4}\frac{2-p}{p-1}}
\\
\\
+ \dfrac{M^{\frac{1}{2}+\frac{1}{4}\frac{2-p}{p-1}}}{\sqrt {t-t_0+T}} 
+ \sup_{Q_{R, T}} \left\{\left[v^{\frac{1}{4}\frac{p+2}{p-1}} |\Sigma_x(t,x,v)|\right]^{1/3}\right\} 
\\
\\
+ \sup_{Q_{R, T}}\left\{v^\frac{p}{4(p-1)}\left[2 \Sigma_v(t,x,v)+\dfrac{ (2-p)\Sigma(t,x,v)}{2(p-1)v} \right]_+^{1/2}\right\}
\end{array}
\right\}.
\end{align}
\end{remark}

\begin{remark} \label{Sigma-N-u-v-remark}
We can alternatively rewrite \eqref{EQ-eq-2.4} as $\Sigma(t,x,v) = p u^{p-2} {\mathscr N}(t,x,u)$ where $u$ and $v$ 
are related via $v =pu^{p-1}/(p-1)$. In passing we point out that the change of variables taking from $u$ to $v$ is 
often called the {\it pressure} transform in literature. 
\end{remark}

Subject to the bounds in Theorem $\ref{thm1}$ being global in space by passing to the limit $R \to \infty$ we have 
the following global (in space) estimate.

\begin{theorem} \label{thm1-global}
Let $(\mathscr M, g, d\mu)$ be a complete smooth metric measure space with $d\mu=e^{-f} dv_g$. 
Assume the metric and potential are time dependent, of class $\mathscr{C}^2$ and 
that for suitable constants $k, h \ge 0$ and $m \ge n$ satisfy 
${\mathscr Ric}_f^m (g) \ge -(m-1)k g$, $\partial_t g \ge -2 h g$ 
on $\mathscr M \times [t_0-T, t_0]$. Let $u$ be a positive solution 
to \eqref{eq11} with $1 < p <1+ 1/(\sqrt{2m} +1)$ and $v =pu^{p-1}/(p-1)$ and 
$M= \sup v$. Then there exists $C=C(p,\beta,m)>0$ such 
that for every $x \in \mathscr M$ and $t_0-T<t \le t_0$ we have 
\begin{align}\label{eq-2.1-global}
\frac{|\nabla v|}{v^{\beta/2}}(x,t) 
\le C \left \{ \begin {array}{ll}
\sqrt k M^{1-\beta/2} + \left[ \sqrt h + \dfrac{1}{\sqrt {t-t_0+T}} \right] M^{(1-\beta)/2}
\\
\\
+ \sup_{\mathscr M \times [t_0-T, t_0]} \left\{\left[\dfrac{|\Sigma_x(t,x,v)|}{v^{(3\beta-2)/2}}\right]^{1/3}\right\} 
\\
\\
+ \sup_{\mathscr M \times [t_0-T, t_0]}\left\{v^{(1-\beta)/2} \left[2 \Sigma_v(t,x,v)-\dfrac{ \beta \Sigma(t,x,v)}{v} \right]_+^{1/2}\right\}
\end{array}
\right\}.
\end{align}
Here $\beta \in (\beta_1, \beta_2)$ where $\beta_1<\beta_2< 0$
are the roots of $\beta^2 + (2-p)/(p -1)\beta +m/2 =0$ and $\Sigma(t,x,v)$ is as in \eqref{EQ-eq-2.4}. 
\end{theorem}

The special case of time independent metrics and potentials (the static case $\partial_t g \equiv 0$ 
and $\partial_t f \equiv 0$) is of enough significance to be formulated as a separate corollary. Here 
we describe the local version. The global version follows by passing to the limit $R\to\infty$.

\begin{theorem} \label{thm1-static}
Let $(\mathscr M, g, d\mu)$ be a complete smooth metric measure space with $d\mu=e^{-f} dv_g$ and assume 
${\mathscr Ric}_f^m (g) \ge -(m-1)k g$ in ${\mathscr B}_R$ for some $k \ge 0$, $m \ge n$ 
and $R>0$. Let $u$ be a positive solution 
to \eqref{eq11} with $1 < p <1+ 1/(\sqrt{2m} +1)$ and $v =pu^{p-1}/(p-1)$ and 
$M= \sup_{Q_{R,T}} v$. Then there exists $C=C(p,\beta,m)>0$ such 
that for every $(x,t)$ in $Q_{R/2,T}$ with $t>t_0-T$ we have 
\begin{align}\label{eq-2.1-static}
\frac{|\nabla v|}{v^{\beta/2}}(x,t) 
\le C \left \{ \begin {array}{ll}
\left[ \dfrac{k^{1/4}}{\sqrt R} + \dfrac{1}{R} + \sqrt k \right] M^{1-\beta/2}
\\
\\
+ \dfrac{M^{(1-\beta)/2}}{\sqrt {t-t_0+T}}  
+ \sup_{Q_{R, T}} \left\{\left[\dfrac{|\Sigma_x(t,x,v)|}{v^{(3\beta-2)/2}}\right]^{1/3}\right\} 
\\
\\
+ \sup_{Q_{R, T}}\left\{v^{(1-\beta)/2} \left[2 \Sigma_v(t,x,v)-\dfrac{ \beta \Sigma(t,x,v)}{v} \right]_+^{1/2}\right\}
\end{array}
\right\}.
\end{align}
Here $\beta \in (\beta_1, \beta_2)$ where $\beta_1<\beta_2< 0$
are the roots of $\beta^2 + (2-p)/(p -1)\beta +m/2 =0$ and $\Sigma(t,x,v)$ is as in \eqref{EQ-eq-2.4}. 
\end{theorem}

We end this section with an application of the estimates above to parabolic Liouville-type theorems. 
Here by an ancient solution $u=u(x,t)$ to \eqref{eq11} we mean a solution defined on $\mathscr M$ for 
all negative times, that is, for all $(x,t)$ with $x \in \mathscr M$ and $-\infty < t <0$.

\begin{theorem} \label{ancient-1}
Let $(\mathscr M, g, d\mu)$ be a complete smooth metric measure space with $d\mu=e^{-f}dv_g$ and 
${\mathscr Ric}^m_f(g) \ge 0$. Assume $[2(2-p)+\beta(p-1)] {\mathscr N}(u) - 2u{\mathscr N}_u(u) \ge 0$ 
for some $\beta \in (\beta_1, \beta_2)$ and all $u>0$ where $1 < p <1+ 1/(\sqrt{2m} +1)$. 
Then any positive ancient solution $u=u(x,t)$ to the nonlinear porous medium equation 
\begin{equation} \label{ancient-equation}
\partial_t u - \Delta_f u^p = {\mathscr N}(u(x,t)), 
\end{equation}
satisfying the growth at infinity 
\begin{equation}
u(x,t) = o \left( [\varrho(x) + \sqrt{|t|}]^{2/[(p-1)(2-\beta)]} \right),
\end{equation} 
must be spatially constant. In particular, if additionally, ${\mathscr N}(u) \ge a$ for some $a>0$ and all $u>0$ 
then \eqref{ancient-equation} admits no such ancient solutions. 
\end{theorem}

\section{Evolution inequalities ${\bf I}$: $\mathscr L_v^p=\partial_t - (p-1)v\Delta_f$ and $w = |\nabla v|^2/v^{\beta}$} 
\label{sec3}

In this section we derive evolution identities and inequalities for $w = |\nabla v|^2/v^{\beta}$ 
where $v$ relates to the positive solution $u$ through the pressure transform $v =pu^{p-1}/(p-1)$. 
Here the evolution operator is $\mathscr L^p_v = \partial_t -(p-1) v \Delta_f$ and $\beta \in {\mathbb R}$ 
is to be specified later. %Note that initially we only assume $p>1$ but later in the course of the proof 
%of the estimates and bounds (see the next section) we need to further restrict the range of $p$. 
\begin{lemma}\label{Lem.2.2}
Let $u$ be a positive solution to \eqref{eq11} with $p>1$ and let $v =pu^{p-1}/(p-1)$ and 
$\Sigma=\Sigma(t,x,v)$ be as in \eqref {EQ-eq-2.4}. Then $v$ satisfies the evolution equation 
\begin{align}\label{eq-2.9}
\mathscr L_v^p [v] = [\partial_t -(p-1) v \Delta_f] v = |\nabla v|^2 + \Sigma(t,x,v). 
%p [(p-1)v/p ]^{\frac{p -2}{p-1}} \mathscr N \left(t,x,[(p-1)v/p]^{\frac{1}{p-1}}\right).
\end{align}
\end{lemma}

\begin {proof}
Using the formulation of $v$ and by directly differentiating $u$ and $u^p$ we have 
\begin{align}
\partial_t u =&~\partial_t [(p-1)v/p]^{1/(p-1)} = (1/p) [(p-1)v/p]^{\frac{2-p}{p-1}} \partial_t v, \nonumber \\ 
\nabla u^p =&~\nabla [(p-1)v/p]^{p/(p-1)} = [(p-1)v/p]^{\frac{1}{p-1}} \nabla v, \nonumber \\
\Delta u^p 
%=&~ [(p-1)v/p]^{\frac{1}{p-1}} \Delta v + (1/p) [(p-1)v/p ]^{\frac{2-p}{p-1}} |\nabla v|^2\nonumber\\
=&~(1/p) [(p-1)v/p]^{\frac{2-p}{p-1}}[(p-1) v \Delta v + |\nabla v|^2].
\end{align}
Therefore by recalling $\Delta_f u^p = \Delta u^p - \langle \nabla f , \nabla u^p \rangle$ and upon substitution it is seen that 
\begin{align}
\Delta_f u^p 
&= (1/p) [(p-1)v/p]^{\frac{2-p}{p-1}}[(p-1) v \Delta v + |\nabla v|^2]
- [(p-1)v/p]^{\frac{1}{p-1}}\langle \nabla f, \nabla v \rangle\nonumber\\
&= (1/p) [(p-1)v/p]^{\frac{2-p}{p-1}}[(p-1) v \Delta_f v + |\nabla v|^2].
\end{align}
Now referring to equation \eqref{eq11} and by substituting the above fragments in the equation we have 
\begin{align}
\square_p u = \partial_t u - \Delta_f u^p &= (1/p) [(p-1)v/p]^{\frac{2-p}{p-1}} (\partial_t v - (p-1)v \Delta_f v 
- |\nabla v|^2) \nonumber \\
&= \mathscr N\left(t,x,[(p-1)v/p]^{\frac{1}{p-1}}\right), 
\end{align}
which upon a rearrangement of terms and using \eqref {EQ-eq-2.4} gives the desired conclusion. 
\end{proof}

\begin{lemma} \label{Lem-2.3}
Let $u$ be a positive solution to \eqref{eq11} with $p>1$ and set 
$w = |\nabla v|^2/v^{\beta}$ where $v =pu^{p-1}/(p-1)$ and $\beta \in {\mathbb R}$. Then $w$ 
satisfies the evolution equation
\begin{align}\label{eq2.3}
\mathscr L_v^p [w] 
=& -2\left[ \frac{1}{2} \partial_t g+(p -1)v{\mathscr Ric}_f^m (g) \right] \frac{( \nabla v, \nabla v)}{v^\beta}\\
&+\frac{2(p -1)}{v^{\beta}} \left[|\nabla v|^2 \Delta_f v - v|\nabla\nabla v|^2
-\frac{v \langle \nabla f, \nabla v \rangle^2}{(m-n)}\right]\nonumber\\
&+2[1+\beta(p -1)] \frac{\langle \nabla v, \nabla |\nabla v|^2 \rangle}{v^\beta}
-\beta \frac{ |\nabla v|^2 \Sigma(t,x,v)}{v^{\beta+1}}\nonumber\\
&-\beta [1+(p -1)(\beta+1)]\frac{|\nabla v|^4}{v^{\beta+1}}
+2 \frac{\langle \nabla v, \nabla\Sigma(t,x,v)\rangle}{v^{\beta}}. \nonumber
\end{align} 
 \end{lemma}
 
\begin{proof}
As $w = |\nabla v|^2/v^{\beta}$ and $v =pu^{p-1}/(p-1)$ satisfies \eqref{eq-2.9}, a direct calculation gives
\begin{align} \label{eq-8.3}
\partial_t w= \partial_t \left[ \frac{|\nabla v|^2}{v^{\beta}}\right] 
=&~2\frac{\langle\nabla v, \nabla \partial_t v \rangle }{v^\beta}
-\frac{[\partial_t g] ( \nabla v, \nabla v)}{v^\beta}
- \beta\frac{v^{\beta-1}|\nabla v|^2 \partial_t v}{v^{2\beta}}\nonumber\\
=&~ \frac{2\langle \nabla v, \nabla [(p -1) v \Delta_f v +|\nabla v|^2 +\Sigma(t,x,v)]\rangle }{v^\beta}
-\frac{[\partial_t g] ( \nabla v, \nabla v)}{v^\beta} \nonumber\\
&-\frac{\beta |\nabla v|^2 [(p -1) v \Delta_f v +|\nabla v|^2 +\Sigma(t,x,v)]}{v^{\beta+1}}\nonumber\\
=& ~ 2(p -1)\frac{ |\nabla v|^2 \Delta_f v}{v^\beta}
+2(p -1) \frac{\langle \nabla v, \nabla \Delta_f v \rangle}{v^{\beta -1}}
+2\frac{\langle \nabla v , \nabla |\nabla v|^2 \rangle}{v^\beta}\nonumber\\
&+ 2 \frac{\langle \nabla v, \nabla \Sigma(t,x,v)\rangle}{v^\beta} 
-\frac{[\partial_t g]( \nabla v, \nabla v)}{v^\beta}
- \beta(p -1) \frac{|\nabla v|^2 \Delta_f v}{v^\beta}\nonumber\\
& - \beta\frac{|\nabla v|^4}{v^{\beta+1}}
- \beta\frac{|\nabla v|^2 \Sigma(t,x,v)}{v^{\beta+1}}.
\end{align}

Furthermore, as $\nabla w = \nabla ( |\nabla v|^2/v^{\beta})= \nabla |\nabla v|^2/v^\beta
- \beta |\nabla v|^2 \nabla v/v^{\beta+1}$, 
by recalling the formulation $\Delta_f w = \Delta w - \langle \nabla f, \nabla w \rangle$ it is seen that 
\begin{align} \label{eq-8.6}
\Delta_f w &= {\text div} [\nabla |\nabla v|^2/v^\beta
- \beta |\nabla v|^2 \nabla v/v^{\beta+1} ]- \langle \nabla f, \nabla |\nabla v|^2/v^\beta
- \beta |\nabla v|^2 \nabla v/v^{\beta+1} \rangle \nonumber\\
&=\frac{\Delta_f |\nabla v|^2}{v^\beta}
- 2 \beta\frac{\langle \nabla v, \nabla |\nabla v|^2 \rangle}{v^{\beta+1}}
-\beta \frac{|\nabla v|^2 \Delta_f v}{v^{\beta+1}}
+\beta(\beta+1) \frac{ |\nabla v|^4}{v^{\beta+2}}.
\end{align} 
Putting \eqref{eq-8.3}-\eqref{eq-8.6} together and rearranging terms results in
\begin{align}\label{eq-8.8}
[\partial_t-(p-1) v\Delta_f] w = &
~2(p -1) \frac{|\nabla v|^2 \Delta_f v}{v^{\beta}}
+2[1+ \beta(p -1)] \frac{\langle \nabla v, \nabla |\nabla v|^2 \rangle}{v^\beta}\nonumber\\
&+2(p -1) \frac{\langle \nabla v, \nabla \Delta_f v \rangle}{v^{\beta -1}}
-(p -1) \frac{\Delta_f |\nabla v|^2}{v^{\beta -1}}\nonumber\\
&-[\beta+(p -1) \beta(\beta+1)]\frac{|\nabla v|^4}{v^{\beta+1}}
-\frac{[\partial_t g] ( \nabla v, \nabla v)}{v^\beta}\nonumber\\
&+2 \frac{\langle \nabla v, \nabla \Sigma(t,x,v) \rangle}{v^{\beta}}
- \beta\frac{ |\nabla v|^2 \Sigma(t,x,v)}{v^{\beta+1}}.
\end{align}
Applying the weighted Bochner-Weitzenb\"ock formula \eqref{Bochner-1} to the right-hand side gives
\begin{align}\label{eq8.10}
[\partial_t-(p-1) v\Delta_f] w = &
-\frac{2(p -1)}{v^{\beta -1}}\left[|\nabla\nabla v|^2
+{\mathscr Ric}_f^m (\nabla v, \nabla v)
+\frac{\langle \nabla f, \nabla v \rangle^2}{(m-n)}\right]\nonumber\\
&+ 2(p -1) \frac{|\nabla v|^2 \Delta_f v}{v^{\beta}}
+ 2[1+\beta(p -1)] \frac{\langle \nabla v, \nabla |\nabla v|^2 \rangle}{v^\beta}\nonumber\\
&- [\beta+(p -1) \beta(\beta+1)]\frac{|\nabla v|^4}{v^{\beta+1}}
- \frac{[\partial_t g] ( \nabla v, \nabla v)}{v^\beta}\nonumber\\
&+ 2 \frac{\langle \nabla v, \nabla\Sigma(t,x,v)\rangle}{v^{\beta}}
- \beta \frac{|\nabla v|^2 \Sigma(t,x,v)}{v^{\beta+1}},
\end{align}
which after rearranging term leads to the desired identity. 
\end{proof}

\begin {lemma} \label{Lem-8.2}
Under the assumptions of Lemma $\ref {Lem-2.3}$, if the metric $g$ and potential $f$ 
satisfy the super flow inequality  
\begin{align} \label{flow-inequality-SPR-PME}
\frac{1}{2} \partial_t g + (p-1) v {\mathscr Ric}_f^m (g) \ge - {\mathsf k} g,  
\end{align}
then $w$ satisfies the evolution inequality 
\begin{align} \label{Eq-8.2}
\mathscr L_v^p [w] \le&~ 
2 {\mathsf k} w + 2 [1+\beta(p -1)] \langle \nabla v , \nabla w \rangle \\
& + ( p -1) \left[ \beta ^2 -\frac{p -2}{p -1} \beta +\frac{m}{2} \right] v^{\beta-1} w^2 \nonumber \\
&+ 2 \frac{\langle \nabla v, \Sigma_x(t,x,v) \rangle}{v^\beta} +\left[2 \Sigma_v(t,x,v)
-\frac{ \beta \Sigma(t,x,v)}{v} \right]w. \nonumber
\end{align}
\end{lemma}

Note that $\mathsf k$ here may be a constant or more generally a function of space and time.

\begin{proof} 
Starting from \eqref{eq2.3} and noting
$\langle \nabla v, \nabla w \rangle = \langle \nabla v, \nabla |\nabla v|^2 \rangle/v^\beta
-\beta |\nabla v|^4/v^{\beta+1}$ we can write upon substitution 
\begin{align}
[\partial_t-(p-1) v\Delta_f] w = 
& -2\left[\frac{1}{2} \partial_t g+(p -1)v{\mathscr Ric}_f^m(g)\right] \frac{( \nabla v, \nabla v)}{v^\beta}\\
&+\frac{2(p -1)}{v^{\beta}} \left[|\nabla v|^2 \Delta_f v - v|\nabla\nabla v|^2
-\frac{v \langle \nabla f, \nabla v \rangle^2}{(m-n)}\right]\nonumber\\
&+2[1+\beta(p -1)] \langle \nabla v, \nabla w \rangle 
-\beta \frac{ |\nabla v|^2 \Sigma(t,x,v)}{v^{\beta+1}}\nonumber\\
&+[(p-1)\beta^2-(p-2)\beta]\frac{|\nabla v|^4}{v^{\beta+1}}
+2 \frac{\langle \nabla v, \nabla\Sigma(t,x,v)\rangle}{v^{\beta}}.\nonumber
\end{align}
Now by writing $\langle \nabla v, \nabla \Sigma(t,x,v) \rangle 
=\langle \nabla v, \Sigma_x(t,x,v) \rangle + \Sigma_v(t,x,v) |\nabla v|^2$ 
and substituting this back in the equation it follows that
\begin{align}\label{E-eq-2.23}
[\partial_t-(p-1) v\Delta_f] w =
& -2\left[\frac{1}{2}\partial_t g+(p -1)v{\mathscr Ric}_f^m(g)\right] \frac{( \nabla v, \nabla v)}{v^\beta}\\
&+\frac{2(p -1)}{v^{\beta}} \left[|\nabla v|^2 \Delta_f v - v|\nabla\nabla v|^2
-\frac{v \langle \nabla f, \nabla v \rangle^2}{(m-n)}\right]\nonumber\\
&+2[1+\beta(p -1)] \langle \nabla v, \nabla w \rangle 
+\left[2 \Sigma_v(t,x,v) -\beta \frac{ \Sigma(t,x,v)}{v}\right] \frac{|\nabla v|^2}{v^{\beta}}\nonumber\\
&+[(p-1)\beta^2-(p-2)\beta]\frac{|\nabla v|^4}{v^{\beta+1}}
+2\frac{\langle \nabla v, \Sigma_x(t,x,v) \rangle}{v^\beta}.\nonumber
\end{align}
Next referring to the second line in the above and by making note of the inequality
\begin {equation} \label{EQ-eq2.24-PME}
|\nabla \nabla v|^2 + \frac{\langle \nabla f, \nabla v \rangle^2}{m-n} 
\ge \frac{(\Delta v)^2}{n} + \frac{\langle \nabla f, \nabla v \rangle ^2}{m-n} 
\ge \frac{(\Delta_f v)^2}{m},  
\end{equation}
we can write 
\begin{align} \label{eq-8.12}
\frac{2(p -1)}{v^{\beta}} &\left[|\nabla v|^2 \Delta_f v
-v|\nabla\nabla v|^2 - \frac{v\langle \nabla f, \nabla v \rangle^2}{m-n}\right] \nonumber\\
&\le 2(p -1) \frac{1}{v^{\beta -1}} \left[ \frac{|\nabla v|^2 \Delta_f v}{v}- \frac{(\Delta_f v)^2}{m}\right] \\
&\le 2(p -1) \frac{1}{v^{\beta -1}} \left[-\left( \frac{\sqrt m}{2} \frac{|\nabla v|^2}{v}  - \frac{\Delta_f v}{\sqrt m} \right)^2
+\frac{m}{4} \frac{|\nabla v|^4}{v^2} \right] 
\le \frac{m(p -1)}{2} \frac{|\nabla v|^4}{v^{\beta+1}}. \nonumber 
\end{align}
Therefore substituting the above inequality in \eqref{E-eq-2.23} results in
\begin{align}\label{eq8.13}
{\mathscr L}_v^p[w] \le &
-2\left[\frac{1}{2}\partial_t g+(p -1)v{\mathscr Ric}_f^m (g)\right]\frac{( \nabla v, \nabla v)}{v^\beta} 
+2[1+ \beta(p -1)] \langle \nabla v, \nabla w \rangle\nonumber\\
&+ \left[\frac{m}{2} (p -1)+ (p-1)\beta^2-(p-2)\beta \right] \frac{|\nabla v|^4}{v^{\beta+1}} \\
&+2\frac{\langle \nabla v, \Sigma_x(t,x,v) \rangle}{v^\beta}
+\left[2 \Sigma_v(t,x,v) -\beta \frac{ \Sigma(t,x,v)}{v}\right]\frac{|\nabla v|^2}{v^{\beta}}. \nonumber 
\end{align}
The conclusion now follows by making use of the super flow inequality \eqref{flow-inequality-SPR-PME}, 
making the substitution $w = |\nabla v|^2/v^{\beta}$ and rearranging terms. 
\end{proof}

\section{Space-time cut-offs, cylindrical localisation and the proof of Theorem \ref{thm1}}
\label{sec4}

We start this section with a construction of smooth space-time cut-off functions that will be used 
in the proof of Theorem \ref{thm1}. First, let us fix a reference point $(x_0,t_0)$ with $x_0 \in \mathscr M$, 
$t_0 \in {\mathbb R}$, and pick $R, T>0$ and then $\tau \in (t_0-T, t_0]$. The following standard lemma grants 
the existence of a smooth {\it profile} function $\bar{\eta}=\bar \eta(\varrho, t)$ with variables $\varrho \ge 0$ 
and $t_0-T \le t \le t_0$ that satisfies a suitable set of properties and bounds.

\begin{lemma} \label{phi lemma} Fix $t_0 \in {\mathbb R}$ and let $R, T>0$. Given $\tau \in (t_0-T, t_0]$ 
there exists a smooth function $\bar{\eta}:[0,\infty) \times [t_0-T, t_0] \to \mathbb{R}$ such that the following 
properties hold:
\begin{enumerate}[label=$(\roman*)$]
\item ${\rm supp} \, \bar{\eta}(\varrho,t) \subset [0,R] \times [t_0-T, t_0]$  and $0 \leq \bar{\eta}(\varrho,t) \leq 1$ in $[0,R] \times [t_0-T, t_0]$,
\item $\bar{\eta}=1$ in $[0,R/2] \times [\tau, t_0]$ and $\partial \bar{\eta}/\partial \varrho =0$ in $[0,R/2] \times [t_0-T, t_0]$, respectively,
\item there exists $c>0$ such that 
\begin{equation}
\frac{|\partial_t \bar{\eta}|}{\sqrt{\bar\eta}} \le \frac{c}{\tau-t_0+T}, 
\end{equation} 
in $[0,\infty)\times[t_0-T,t_0]$ and $\bar{\eta}(\varrho,t_0-T)=0$ for all $\varrho \in [0,\infty)$.
\item $-c_a \bar{\eta}^a/R \le \partial_\varrho \bar{\eta} \le 0$ and $|\partial_{\varrho \varrho} \bar{\eta}| \le c_a \bar{\eta}^a / R^2$ 
hold on $[0, \infty)\times [t_0-T, t_0]$ for every $0<a<1$ and some $c_a>0$.
\end{enumerate}
\end{lemma}

Using this lemma we next introduce a smooth space-time cut-off function $\eta=\eta(x,t)$ by way of composition 
of the profile $\bar{\eta}$ with the geodesic radial variable $\varrho$, namely, 
\begin{equation} \label{cut-off def}
\eta(x,t) = \bar{\eta}(\varrho(x,t), t), \qquad (x, t) \in \mathscr M \times [t_0-T, t_0]. 
\end{equation}

It is clear that $0 \le \eta \le 1$ and ${\rm supp}\,\eta \subset Q_{R,T} = \{(x,t) | d(x,x_0, t) \le R, t_0-T \le t \le t_0\}$ whilst $\eta \equiv 1$ 
on $\{(x,t) | d(x,x_0, t) \le R/2, \tau \le t \le t_0\}$. Moreover, referring to \eqref{cut-off def} a straightforward calculation gives  
\begin{itemize}
\item $\nabla \eta = \partial_\varrho \bar\eta \nabla \varrho$, 
\item $\partial_t \eta = \partial_\varrho \bar\eta \partial_t \varrho + \partial_t \bar\eta$, 
\item $\Delta_f \eta = \partial_{\varrho \varrho} \bar\eta |\nabla \varrho|^2 
+ \partial_\varrho \bar\eta \Delta_f \varrho$. 
\end{itemize}
As a result it is easily seen that  
\begin{align}
\mathscr L_v^p [\eta] &= [\partial_t - (p -1)v \Delta_f] \eta \nonumber \\
&= \partial_\varrho \bar\eta \partial_t \varrho + \partial_t \bar\eta 
- (p-1) v [\partial_{\varrho \varrho} \bar\eta |\nabla \varrho|^2 + \partial_\varrho \bar\eta \Delta_f \varrho] \nonumber \\
&= \partial_\varrho \bar\eta \mathscr L_v^p [\varrho] + \partial_t \bar\eta - (p-1) v \partial_{\varrho \varrho} \bar\eta |\nabla \varrho|^2.
\end{align}
Another auxiliary result on $\mathscr L_v^p$ that will be needed later on is described below.

\begin{lemma}\label{Lem-2.6}
For functions $u=u(x,t)$ and $w=w(x,t)$ of class $\mathscr C^2$ we have 
\begin{align}
u\mathscr L_v^p[uw] = uw\mathscr L_v^p[u]
- 2(p -1)v [\langle\nabla u, \nabla (u w) \rangle-|\nabla u|^2 w]+u^2\mathscr L_v^p [w].
\end{align}
\end{lemma}

\begin{proof}
Using $\partial_t (uw) = w \partial_t u + u \partial_t w$ and 
$\Delta_f (uw) = w \Delta_f u + 2 \langle \nabla u, \nabla w \rangle + u \Delta_f w$ we can write by substitution 
\begin{align}
\mathscr L_v^p [u w] &= [\partial_t-(p -1)v \Delta_f] (u w) \nonumber \\
&= w \partial_t u + u \partial_t w - (p-1) v [w \Delta_f u + 2 \langle \nabla u, \nabla w \rangle + u \Delta_f w] \nonumber\\
&= w [\partial_t-(p -1)v \Delta_f] u-2(p-1)v\langle \nabla u, \nabla w \rangle 
+u[\partial_t-(p -1)v \Delta_f] w.
\end{align}
Now multiplying through by $u$ and using 
$u \langle \nabla u, \nabla w \rangle = \langle \nabla u, \nabla (uw) \rangle - w |\nabla u|^2$ gives
\begin{align}
u\mathscr L_v^p [u w] =&~ uw [\partial_t-(p -1)v \Delta_f] u 
- 2(p -1)v [\langle\nabla u ,\nabla (u w) \rangle- |\nabla u|^2 w] \nonumber\\
&+u^2[\partial_t-(p -1)v \Delta_f] w,
\end{align}
which is the desired identity.
\end{proof}

\qquad \\
{\bf Proof of Theorem \ref{thm1}.} We begin by considering the localised function $\eta w$ where $\eta$ 
is the cut-off function from \eqref{cut-off def}. An application of Lemma \ref{Lem-2.6} then gives 
\begin{align}
\mathscr L_v^p[\eta w] = w\mathscr L_v^p[\eta]
- 2(p -1)v [\langle\nabla \eta ,\nabla (\eta w) \rangle- |\nabla \eta|^2 w]/\eta
+\eta \mathscr L_v^p [w].
\end{align}

We now wish to utilise the bound on $\mathscr L_v^p[w]$ from \eqref{Eq-8.2} in Lemma \ref{Lem-8.2}. Towards this end 
invoking the assumptions in Theorem \ref{thm1} we have the super flow inequality \eqref{flow-inequality-SPR-PME} 
with ${\mathsf k} = (p-1)(m-1) k v + h$ [we can also take ${\mathsf k} = (p-1)(m-1) k M + h$ as $0<v \le M$]. Therefore 
substituting this together with 
$\eta \langle \nabla v, \nabla w \rangle = \langle \nabla v, \nabla (\eta w) \rangle - w \langle \nabla v, \nabla \eta \rangle$ 
then leads to
\begin{align}\label{eq8.18}
\mathscr L_v^p [\eta w] \le &~
w \mathscr L_v^p [\eta]
-2(p -1) v \left\langle \frac{\nabla \eta}{\eta} , \nabla (\eta w) \right\rangle
+ 2(p -1)v \frac{|\nabla \eta|^2}{\eta} w \nonumber\\
&+ 2[(p - 1) (m-1)k v + h] \eta w
+2 [1+\beta(p -1)] \langle \nabla v , \nabla (\eta w) \rangle\nonumber\\
&-2 [1+\beta(p -1)] w\langle \nabla v , \nabla \eta \rangle
+ ( p -1) \left[ \beta ^2 -\frac{p -2}{p -1} \beta +\frac{m}{2} \right] v^{\beta-1} \eta w^2\nonumber\\
&+ 2\eta \frac{\langle \nabla v, \Sigma_x(t,x,v) \rangle}{v^\beta} + \left[2 \Sigma_v(t,x,v)
-\frac{ \beta \Sigma(t,x,v)}{v} \right] \eta w.
\end{align}  
Assume now that the localised function $\eta w$ is maximised at the point $(x_1, t_1)$ in the compact set 
$\{d(x,x_0, t) \le R, t_0-T \le t \le \tau\} \subset \mathscr M \times [t_0-T, t_0]$. Additionally, we can assume that $(\eta w)(x_1, t_1) >0$ 
as otherwise the conclusion of the theorem is true with $w(x, \tau) \le 0$ for all $d(x, x_0, \tau) \le R/2$. In particular 
$t_1>t_0-T$ and at the point $(x_1,t_1)$  we have $\Delta_f(\eta w) \le 0$, $\partial_t (\eta w) \ge 0$ and 
$\nabla(\eta w) =0$. Thus at the point $(x_1,t_1)$ we also have $\mathscr L_v^p [\eta w] \ge 0$. 
From \eqref{eq8.18} upon rearranging of terms it therefore follows that  
\begin{align}\label{E-eq2.29}
- ( p -1) \left[ \beta ^2 -\frac{p -2}{p -1} \beta 
+\frac{m}{2} \right] v^{\beta-1} \eta w^2 \le
&~ 2[(p - 1)(m-1) k v + h] \eta w\nonumber\\
&- 2 [1+\beta(p -1)] w\langle \nabla v , \nabla \eta \rangle \nonumber\\
&+ 2(p -1)v \frac{|\nabla \eta|^2}{\eta} w
+ w \mathscr L_v^p [\eta]\nonumber\\
&+\left[2 \Sigma_v(t,x,v)-\frac{ \beta \Sigma(t,x,v)}{v} \right] \eta w\nonumber\\
&+ 2\eta \frac{\langle \nabla v, \Sigma_x(t,x,v) \rangle}{v^\beta}.
\end{align}
Now it is easily seen that the quadratic expression $\beta ^2 -(p -2)/(p -1) \beta +m/2$ has the discriminant $[(p -2)/(p -1)]^2 -2m$ 
which is positive only when $p \in (1,1+1/[\sqrt{2m}+1])$. Thus denoting the roots of this quadratic expression by 
$\beta_1<\beta_2< 0$, it is evident upon choosing $\beta \in (\beta_1, \beta_2)$ that this quadratic expression will be negative. 
With this choice of $\beta$ we now set $-(p -1)[\beta ^2 -(p -2)/(p -1) \beta +m/2] =2/ \gamma$ where $\gamma > 0$.
Therefore returning to \eqref{E-eq2.29} and multiplying through by $\gamma v^{1-\beta}$ we can write 
\begin{align}\label{eq-8.20}
2 \eta w^2 \le&~ 
2 \gamma [(p - 1)(m-1) k v^{2-\beta} + h v^{1-\beta}]\eta w\nonumber\\
&- 2 \gamma [1+\beta(p -1)] v^{1-\beta}w\langle \nabla v , \nabla \eta \rangle\nonumber\\
&+ 2 \gamma (p -1) v^{2-\beta} \frac{|\nabla \eta|^2}{\eta} w
+ \gamma v^{1-\beta} w \mathscr L_v^p [\eta]\nonumber\\
&+\left[2 \Sigma_v(t,x,v)-\frac{ \beta \Sigma(t,x,v)}{v} \right] \gamma v^{1-\beta}\eta w\nonumber\\
&+ 2 \gamma v^{1-\beta} \eta \frac{\langle \nabla v, \Sigma_x(t,x,v) \rangle}{v^\beta}. 
\end{align}

Now we proceed onto bounding the expression on the right-hand side of \eqref{eq-8.20}.
Towards this end, starting from the first line, upon recalling $0\le \eta \le 1$ we can write 
\begin{align}
2 \gamma (p - 1)(m-1) k v^{2-\beta} \eta w &\le
2 \gamma (p - 1)(m-1) k v^{2-\beta} \eta^{1/2} w\\
&\le \frac{1}{7} \eta w^2 + C \gamma ^2 (p-1)^2 (m-1)^2 k^2(\sup_{Q_{R, T}} v)^{4-2\beta},\nonumber
\end{align} 
and 
\begin{align}
2 \gamma h v^{1-\beta} \eta w &\le
2 \gamma h v^{1-\beta} \eta^{1/2} w 
\le \frac{1}{7} \eta w^2 + C \gamma ^2 h^2 (\sup_{Q_{R, T}} v)^{2-2\beta}.
\end{align} 
Note that in view of $\beta<0$ we have $2- \beta>0$ and $1- \beta>0$. 
Moving to the second line, by recalling $w = |\nabla v|^2/v^{\beta}$ and using 
$(iv)$ in Lemma \ref{phi lemma}, we have 
\begin{align}\label{bound1}
-2 \gamma [1+\beta(p -1)] v^{1-\beta}w \langle \nabla v , \nabla \eta \rangle 
& \le 2 \gamma |1+\beta(p -1)| v^{1-\beta}w |\nabla v| |\nabla \eta|\nonumber\\
& \le 2 \gamma |1+\beta(p -1)| v^{1-\beta/2} \eta^{3/4} w^{3/2} \frac{|\nabla \eta|}{\eta^{3/4}}\nonumber\\
& \le \frac{1}{7} (\eta ^{3/4} w^{3/2})^{4/3} + C \gamma^4 [1+\beta(p -1)]^4 \left[\frac{|\nabla \eta|}{\eta^{3/4}}\right]^4 v^{4-2\beta}\nonumber\\
& \le \frac{1}{7} \eta w^2 + C \gamma^4\frac{[1+\beta(p -1)]^4}{R^4} (\sup_{Q_{R, T}} v)^{4-2\beta}.
\end{align}
In much the same way for the next term by taking advantage of $(iv)$ in Lemma \ref{phi lemma} we can write
 \begin{align}
 2 \gamma (p -1) v^{2-\beta} \frac{|\nabla \eta|^2}{\eta} w 
 &= 2\gamma (p -1) v^{2-\beta} \frac{|\nabla \eta|^2}{\eta^{3/2}} \eta^{1/2} w \nonumber\\
 & \le \frac{1}{7} \eta w^2 +C \gamma^2 (p-1)^2 \left[\frac{|\nabla \eta|}{\eta ^{3/4}} \right]^4 v^{4-2\beta}\nonumber\\
 & \le \frac{1}{7} \eta w^2 + C \gamma^2\frac{(p -1)^2}{R^4}(\sup_{Q_{R, T}} v)^{4-2\beta}.
 \end{align}

We now come to the term 
$\gamma v^{1-\beta} w \mathscr L_v^p [\eta] = \gamma v^{1-\beta} w[\partial_t - (p -1)v \Delta_f] \eta$. 
To handle this term we consider the space and time derivatives separately and provide a lower bound for each.
\begin{itemize}
\item Bounding ${\bf I}=\gamma (p-1) v^{2-\beta} w (-\Delta_f \eta)$: 
By virtue of ${\mathscr Ric}_f^m(g) \geq - (m-1)k g$ we have
$\Delta_f \varrho \le(m-1) \sqrt {k}\coth (\sqrt{k} \varrho)$
and so from \eqref{cut-off def}, $\eta$ being radial and
$ \bar\eta_\varrho \le 0$, it follows that
\begin{align}\label{eq2.26}
\Delta_f \eta = \bar\eta_{\varrho \varrho}|\nabla \varrho|^2 + \bar\eta_\varrho \Delta_f \varrho
\ge \bar\eta_{\varrho \varrho} + \bar\eta_\varrho (m-1)\sqrt{k} \coth (\sqrt{k} \varrho).
\end{align} 
Now using the bound $\sqrt{k} \coth (\sqrt{k} \varrho)\le \sqrt{k} \coth(\sqrt {k} R/2) \le (2+\sqrt {k} R)/R $
for $0 \le \varrho \le R/2$ and noting that $ \bar\eta_\varrho = 0$ for  $0 \le \varrho \le R/2$, it follows that 
\begin{align}\label{eq2.27}
-\Delta_f \eta& \le - [ \bar\eta_{\varrho \varrho} + \bar\eta_\varrho (m-1)\sqrt {k} \coth(\sqrt {k} R/2)]\nonumber\\
&\le |\bar\eta_{\varrho \varrho}| + (m-1)\left(\frac{2}{R} +\sqrt {k}\right) |\bar\eta_\varrho|.
\end{align} 
Therefore using Young's inequality and $(iv)$ in Lemma \ref{phi lemma} we can write 
\begin{align}\label{eq2.38}
&\gamma (p -1) v^{2-\beta} w (-\Delta_f \eta)\\
& \le \gamma (p -1)  v^{2-\beta} w \left[ |\bar\eta_{\varrho \varrho}| 
+(m-1)\left(\frac{2}{R} +\sqrt {k} \right) |\bar\eta_\varrho| \right]\nonumber\\
& \le\gamma (p -1) v^{2-\beta}  \sqrt \eta w 
 \left[\frac{|\bar\eta_{\varrho \varrho}|}{\sqrt \eta}
+(m-1)\left(\frac{2}{R} +\sqrt {k} \right) \frac{|\bar\eta_\varrho|}{\sqrt \eta} \right]\nonumber\\
& \le \frac{1}{14} \eta w^2 + C\gamma^2 (p -1)^2
\left[\left(\frac{|\bar\eta_{\varrho \varrho}|}{\sqrt \eta}\right)^2
+(m-1)^2\left(\frac{1}{R^2} +k \right) \left(\frac{|\bar\eta_\varrho|}{\sqrt \eta}\right)^2 \right] v^{4-2\beta}\nonumber\\
& \le \frac{1}{14} \eta w^2 + C \gamma^2(p-1)^2(m -1)^2
\left[\frac{1+ k R^2}{R^4}\right] (\sup_{Q_{R, T}} v)^{4-2\beta}. \nonumber 
\end{align}

\item Bounding ${\bf II}=\gamma v^{1-\beta}w \partial_t \eta$:
We use Lemma \ref{phi lemma} and $\partial_t \varrho (x,t) \ge - h R$ in $Q_{R,T}$. 
Let us first justify the latter inequality. 
Fix $x$ and $t$ such that $d(x, x_0, t)<R$. Let $X(x_0,x)$ be the set 
of all minimal geodesics $\lambda=\lambda(s): [0, 1] \to \mathscr M$ with respect to $g(t)$ connecting the reference 
point $x_0 =\lambda(0)$ to $x=\lambda(1)$ and let $\Omega(x_0,x)$ be the set of all ${\mathscr C}^1$ 
curves connecting $x_0$ to $x$. By using Lemma {\rm B.40} p.~531 in \cite{Chow} we can write
\begin{align} 
\partial_t \varrho (x,t) &= \frac{\partial}{\partial t} d (x, x_0; t) 
= \frac{\partial}{\partial t} \left\{ \inf_{\omega \in \Omega(x_0, x)} \int_0^1 |\omega'(s)|_{g(t)} \, ds \right\} \nonumber \\
&= \frac{\partial}{\partial t} \left\{ \inf_{\omega \in \Omega(x_0, x)} \int_0^1 \sqrt{[g(t)] (\omega'(s), \omega'(s))} \, ds \right\} \nonumber \\
&= \inf_{\lambda \in X(x_0,x)} \int_0^1 \frac{[\partial_t g(t)](\lambda'(s), \lambda'(s))}{2  \sqrt{[g(t)] (\lambda'(s), \lambda'(s))}} \, ds \nonumber\\
&=  \inf_{\lambda \in X(x_0,x)} \int_0^1 \frac{[\partial_t g(t)](\lambda'(s), \lambda'(s))}{2|\lambda'(s)|_{g(t)}} \, ds. 
\end{align}
Hence making note of the lower bound $\partial_t g \ge -2h g$ in $Q_{R,T}$ with $h \ge 0$ gives  
\begin{align} \label{rt-bound}
\partial_t \varrho (x,t) 
&\ge \inf_{\lambda  \in X(x_0,x)} \int_0^1 - h |\lambda' (s)|_{g(t)} \, ds 
= \inf_{\lambda \in X(x_0,x)} \left[ - h \underbrace{\int_0^1 |\lambda' (s)|_{g(t)} \, ds}_{=\varrho(x,t)} \right] \nonumber \\
&\ge - h \varrho(x, t) \ge - h R,
\end{align}
which is the desired inequality. 
Next, upon noting $\partial_t \eta = \bar{\eta}_t + \bar \eta_\varrho \partial_t \varrho$ and $(iii)$ in Lemma \ref{phi lemma} we have 
\begin{align}\label{eq2.30}
\partial_t \eta &= \bar{\eta}_t + \bar \eta_\varrho \partial_t \varrho 
\le |\bar \eta_t| - h R \bar \eta_\varrho \nonumber \\
&\le |\bar \eta_t| + h R |\bar \eta_\varrho| \le C \left[ \frac{1}{\tau-t_0+T} + h \right] \sqrt{\eta}. \nonumber 
\end{align} 
Hence we can write 
\begin{align}
\gamma v^{1-\beta}w \partial_t \eta &= \gamma \sqrt \eta v^{1-\beta} w \frac{\partial_t \eta}{\sqrt \eta}
\le C \gamma \sqrt \eta v^{1-\beta} w \left[\frac{1}{\tau-t_0+T}+h \right]\nonumber\\
&\le \frac{1}{14} \eta w^2 + C \gamma ^2 \left[\frac{1}{(\tau-t_0+T)^2} + h^2\right](\sup_{Q_{R, T}} v)^{2-2\beta}.
\end{align}
\end{itemize}
By putting together the bounds on the fragments ${\bf I}$ and ${\bf II}$ from above we arrive at
\begin{align}\label{EQ-eq-2.46-PME}
\gamma v^{1-\beta} w \mathscr L_v^p [\eta] 
=&~\overbrace{\gamma (p-1) v^{2-\beta} w (-\Delta_f) \eta}^{{\bf I}} + \overbrace{\gamma v^{1-\beta} w\partial_t \eta}^{{\bf II}} \nonumber \\
\le&~\frac{1}{7} \eta w^2
+ C \gamma ^2 \left[\frac{1}{(\tau-t_0+T)^2} + h^2\right](\sup_{Q_{R, T}} v)^{2-2\beta}\nonumber\\
&+ C \gamma^2(p -1)^2 (m-1)^2\left[\frac{1+ k R^2}{R^4}\right] (\sup_{Q_{R, T}} v)^{4-2\beta}.
\end{align}

Finally for the two remaining terms involving the nonlinearity $\Sigma(t,x,v)$, first upon recalling $0\le \eta \le 1$ 
and utilising Young's inequality we have
\begin{align}\label{eq-2.33}
\gamma v^{1-\beta}\eta w \bigg[2 \Sigma_v(t,x,v)&-\frac{ \beta \Sigma(t,x,v)}{v} \bigg] \le 
\gamma v^{1-\beta} \sqrt \eta w \left[2 \Sigma_v(t,x,v)-\frac{ \beta \Sigma(t,x,v)}{v} \right]_+ \nonumber\\
& \le \frac{1}{7} \eta w^2 + C \gamma^2 v^{2(1-\beta)} 
\left[2 \Sigma_v(t,x,v)-\frac{ \beta \Sigma(t,x,v)}{v} \right]_+^2 \\
& \le \frac{1}{7} \eta w^2 + C \gamma^2 \sup_{Q_{R, T}} 
\left\{v^{2(1-\beta)}\left[2 \Sigma_v(t,x,v)-\frac{ \beta \Sigma(t,x,v)}{v} \right]_+^2\right\}. \nonumber 
\end{align}
Likewise, by using the Cauchy-Schwarz and Young inequalities respectively we can bound the last term in \eqref {eq-8.20} as 
\begin{align}
2 \gamma v^{1-\beta} \eta \frac{\langle \nabla v, \Sigma_x(t,x,v) \rangle}{v^\beta} 
\le &~2 \gamma v^{1-\beta} \eta |\nabla v| \frac{|\Sigma_x(t,x,v)|}{v^\beta} \nonumber \\
\le &~ 2 \gamma \eta^{1/4} w^{1/2} v^{1- \beta/2} \frac{|\Sigma_x(t,x,v)|}{v^\beta} \nonumber\\
\le &~\frac{1}{7} \eta w^2 + C \gamma^{4/3} 
\left[ \frac{|\Sigma_x(t,x,v)|}{v^{(3\beta-2)/2}} \right]^{4/3} \nonumber\\
\le &~ \frac{1}{7} \eta w^2+C \gamma^{4/3} \sup_{Q_{R, T}} 
\left\{\left[ \frac{|\Sigma_x(t,x,v)|}{v^{(3\beta-2)/2}} \right]^{4/3}\right\}. 
\end{align}

Completing the estimate of the individual terms on the right-hand side of \eqref{eq-8.20}, by inserting these estimates
back into the inequality and after basic considerations and adjusting constants, we arrive at the following bound 
at the space-time point $(x_1, t_1)$: 
\begin{align}\label{EQ-eq-2.49}
\eta w^2 
\le&~C\gamma^2 
\left \{ \begin {array}{ll} 
(p -1)^2 (m-1)^2\left[\dfrac{1+ k R^2}{R^4}\right]
+\dfrac{(p -1)^2}{R^4}
\\
\\
+\dfrac{\gamma^2 [1+\beta(p -1)]^4}{R^4}
+(p-1)^2 (m-1)^2 k^2
\end{array}
\right\}
(\sup_{Q_{R, T}} v)^{4-2\beta} \nonumber \\
&+C 
\left \{ \begin {array}{ll} 
\gamma^{2} \left[2h^2+\dfrac{1}{(\tau-t_0+T)^2}\right](\sup_{Q_{R, T}} v)^{2-2\beta}
\\
\\
+\gamma^{4/3}\sup_{Q_{R, T}} \left\{\left[ \dfrac{|\Sigma_x(t,x,v)|}{v^{(3\beta-2)/2}}\right]^{4/3}\right\} 
\\
\\
+ \gamma^{2} \sup_{Q_{R, T}} 
\left\{v^{2(1-\beta)}\left[2 \Sigma_v(t,x,v)-\dfrac{ \beta \Sigma(t,x,v)}{v} \right]_+^2\right\} 
\end{array}
\right\}. 
\end{align}

Recalling now that $M=\sup_{Q_{R,T}} v$, and as and a result of $\beta \in (\beta_1, \beta_2)$ 
that $2-2\beta \ge 2$ and $4-2\beta \ge 4$, it follows from \eqref{EQ-eq-2.49} upon absorbing 
$\gamma$ into $C$ and adjusting the constant if necessary that 
\begin{align}
\eta w^2 
\le C(\gamma) & \left \{ \left( \begin {array}{ll}  
(p -1)^2 (m-1)^2\left[\dfrac{1+ k R^2}{R^4}\right] 
+ \dfrac{(p -1)^2}{R^4}
\\
\\
+\dfrac{ [1+\beta(p -1)]^4}{R^4} +(p-1)^2 (m-1)^2 k^2
\end{array}
\right) \right. M^{4-2\beta} \nonumber \\
& +  
\left. \left( \begin {array}{ll} 
\left[ 2h^2 + \dfrac{1}{(\tau-t_0+T)^2} \right] M^{2-2\beta}
%\\
%\\
+ \sup_{Q_{R, T}} \left\{\left[\dfrac{|\Sigma_x(t,x,v)|}{v^{(3\beta-2)/2}}\right]^{4/3}\right\} 
\\
\\
+\sup_{Q_{R, T}} \left\{v^{2(1-\beta)}\left[2 \Sigma_v(t,x,v)-\dfrac{ \beta \Sigma(t,x,v)}{v} \right]_+^2\right\}
\end{array}
\right) \right\}.
\end{align}

Recalling the maximality of $\eta w$ at $(x_1, t_1)$ along with $\eta \equiv 1$ when $d(x, x_0,t) \le R/2$ and $\tau \le t \le t_0$, 
it follows that $w^2 (x, \tau) = (\eta^2 w^2)(x, \tau) \le (\eta^2 w^2)(x_1,t_1) \leq (\eta w^2)(x_1,t_1)$. 
Hence noting $w = |\nabla v|^2/v^{\beta}$ and 
absorbing $\beta$, $p$ and $m$ into $C(\gamma)$ gives
\begin{align}
\frac{|\nabla v|}{v^{\beta/2}} (x,\tau)
\le C(p, \beta, m) \left \{ \begin {array}{ll}
\sqrt h M^{(1-\beta)/2} + \left[ \dfrac{k^{1/4}}{\sqrt {R}} +\dfrac{1}{R} +\sqrt k \right] M^{1-\beta/2}
\\
\\
+ \dfrac{M^{(1-\beta)/2}}{\sqrt {\tau-t_0+T}} 
+ \sup_{Q_{R, T}} \left\{\left[\dfrac{|\Sigma_x(t,x,v)|}{v^{(3\beta-2)/2}}\right]^{1/3}\right\} 
\\
\\
+ \sup_{Q_{R, T}}\left\{v^{(1-\beta)/2}\left[2 \Sigma_v(t,x,v)-\dfrac{ \beta \Sigma(t,x,v)}{v} \right]_+^{1/2}\right\}
\end{array}
\right\}.
\end{align}
The arbitrariness of $\tau$ in the interval $t_0-T<\tau  \le t_0$ now gives the final conclusion.
\hfill $\square$

\section{Ancient solutions to $\partial_t u - \Delta_f u^p = {\mathscr N}(u)$ ${\bf I}$: Proof of Theorem \ref{ancient-1}}
\label{sec5}

Let us now present the proof of Theorem \ref{ancient-1}. Towards this end fix a space-time point $(x_0,t_0)$. 
Then with the choice of base point $(x_0,t_0)$ and $R>0$, $T=R^2$ it follows from the growth assumption 
$u(x,t) = o([\varrho(x) + \sqrt{|t|}]^{2/[(p-1)(2-\beta)]})$ that 
\begin{equation} \label{M-growth-ancient-1}
M = \sup_{Q_{R,T}} v = [p/(p-1)] \sup_{Q_{R,T}} u^{p-1} = o(R^{1/(1-\beta/2)}).
\end{equation} 
Next turning to the local estimate \eqref{eq-2.1-static} in Theorem \ref{thm1-static} (with $t=t_0$, $k=0$ and $T=R^2$) we can write 
\begin{align} \label{vx0t0-ancient-1}
\frac{|\nabla v|}{v^{\beta/2}}(x_0,t_0) 
&\le C \left \{ \begin {array}{ll}
\left[ \dfrac{k^{1/4}}{\sqrt {R}} +\dfrac{1}{R} +\sqrt k \right] M^{1-\beta/2} 
+ \sup_{Q_{R, T}} \left\{\left[\dfrac{|\Sigma_x(t,x,v)|}{v^{(3\beta-2)/2}}\right]^{1/3}\right\} 
\\
\\
+ \dfrac{M^{(1-\beta)/2}}{\sqrt T} 
+ \sup_{Q_{R, T}}\left\{v^{(1-\beta)/2}\left[2 \Sigma_v(t,x,v)-\dfrac{ \beta \Sigma(t,x,v)}{v} \right]_+^{1/2}\right\}
\end{array}
\right\} \nonumber \\
&\le C \left \{ \begin {array}{ll}
\dfrac{M^{1-\beta/2}}{R} + \dfrac{M^{(1-\beta)/2}}{\sqrt T} + \sup_{Q_{R, T}} \left\{\left[\dfrac{|\Sigma_x(t,x,v)|}{v^{(3\beta-2)/2}}\right]^{1/3}\right\} 
\\
\\
+ \sup_{Q_{R, T}}\left\{v^{(1-\beta)/2}\left[2 \Sigma_v(t,x,v)-\dfrac{ \beta \Sigma(t,x,v)}{v} \right]_+^{1/2}\right\}
\end{array}
\right\}.
\end{align}

According to Remark \ref{Sigma-N-u-v-remark}, $\Sigma(t,x,v) = p u^{p-2} {\mathscr N}(t,x,u)$ 
where $u=[(p-1)v/p]^{1/(p-1)}$. Hence, in view of ${\mathscr N}={\mathscr N}(u)$, here, 
$\Sigma=\Sigma(v)$ and so in particular $\Sigma_x \equiv 0$. 
Moreover, a basic calculation based on writing $\Sigma_v=\Sigma_u \partial_v u$ where 
$\Sigma_u = pu^{p-2}[(p-2){\mathscr N}/u+{\mathscr N}_u]$ and $\partial_v u = u^{2-p}/p$ gives 
$\Sigma_v=(p-2){\mathscr N}/u+{\mathscr N}_u$. Therefore by putting the latter together, 
\begin{align}
2 \Sigma_v (v) - \beta \frac{\Sigma(v)}{v} 
&= 2 \left[ (p-2) \frac{{\mathscr N}(u)}{u} + {\mathscr N}_u(u) \right] - \beta (p-1) \frac{{\mathscr N}(u)}{u} \nonumber \\ 
&= [2(p-2)-\beta(p-1)] \frac{{\mathscr N}(u)}{u} + 2 {\mathscr N}_u(u) \le 0,
\end{align}
where the last inequality follows from the assumptions on ${\mathscr N}$ in the theorem. Thus by virtue of 
\eqref{M-growth-ancient-1} and the above calculation, we conclude from \eqref{vx0t0-ancient-1} that  
\begin{align} 
\frac{|\nabla v|}{v^{\beta/2}}(x_0,t_0) \le C\left[ \dfrac{M^{1-\beta/2}}{R} + \dfrac{M^{(1-\beta)/2}}{\sqrt T} \right] 
\le \frac{o(R)}{R} + \frac{o(R^{(1-\beta)/(2-\beta)})}{R}.
\end{align}
Now passing to the limit $R \nearrow \infty$ it follows that $|\nabla v|(x_0,t_0)=0$. The arbitrariness of $(x_0,t_0)$ 
implies $|\nabla v| \equiv 0$ and so $v$ and subsequently $u$ are spatially constant. Hence we have $u=u(t)$. From 
equation \eqref{ancient-equation} satisfied by $u$ it then follows that $du/dt = {\mathscr N}(u)$. Assuming now 
that ${\mathscr N}(u) \ge a >0$ for all $u>0$ it follows by integrating the ODE that $u(t) \le u(0) + at$ for all 
$t<0$. This however clashes with $u(t)>0$ as $t \searrow -\infty$ and so the conclusion is reached. \hfill $\square$

\section {Another Hamilton-Souplet-Zhang estimate for \eqref{eq11}: $1< p < 1+1/\sqrt{m-1}$}
\label{sec6}

In this section we present the second set of gradient estimates for the positive solutions to equation \eqref{eq11}. 
The estimates as will be seen below are valid in the exponent range $1<p<1+1/\sqrt{m-1}$. We use a different 
set of ideas and techniques to prove these estimates, a result of which is that we can improve the key evolution 
inequalities used in the proof, and subsequently arrive at a larger exponent range for the validity 
of the estimates. (Compare with $1<p<1/(\sqrt{2m} + 1)$ in Theorem \ref{thm1}.) 
\begin{theorem} \label{thm1-EQ-PME}
Let $(\mathscr M, g, d\mu)$ be a complete smooth metric measure space with $d\mu=e^{-f} dv_g$. 
Suppose that the metric and potential are time dependent, of class $\mathscr{C}^2$ and that for suitable constants 
$k, h \ge 0$ and $m\ge n$ satisfy the bounds ${\mathscr Ric}_f^m (g) \ge -(m-1)k g$ and 
$\partial_t g \ge -2 h g$ in the space-time cylinder $Q_{R,T}$ with $R, T >0$. Let $u$ be a positive solution 
to \eqref{eq11} with $1 < p < 1+1/\sqrt{m-1}$ and $v =pu^{p-1}/(p-1)$ and $M=\sup_{Q_{R,T}} v$. 
Then there exists $C=C(p,m)>0$ such that for every $(x,t)$ in $Q_{R/2,T}$ with $t>t_0-T$ we have 
\begin{align}\label{eq-2.1-EQ}
v^{\frac{1}{2(p-1)}}|\nabla v| (x,t) 
\le C \left \{ \begin {array}{ll}
\sqrt h M^\frac{p}{2(p-1)} + \left[ \dfrac{k^{1/4}}{\sqrt R} + \dfrac{1}{R} + \sqrt k \right] M^{1+\frac{1}{2(p-1)}}
\\
\\
+ \sup_{Q_{R, T}}\left\{v^{p/[2(p-1)]} \left[2 \Sigma_v(t,x,v)+\dfrac{\Sigma(t,x,v)}{(p-1)v} \right]_+^{1/2}\right\}
\\
\\
+ \dfrac{M^\frac{p}{2(p-1)}}{\sqrt {t-t_0+T}}  
+\sup_{Q_{R, T}} \left\{ \left[ v^{(2p+1)/[2(p-1)]} |\Sigma_x(t,x,v)| \right]^{1/3} \right\} 
\end{array}
\right\}.
\end{align}
\end{theorem}

Subject to the bounds in Theorem $\ref{thm1-EQ-PME}$ being global in space by passing to the limit $R \to \infty$ we have 
the following global (in space) estimate.

\begin{theorem} \label{thm1-EQ-PME-global}
Let $(\mathscr M, g, d\mu)$ be a complete smooth metric measure space with $d\mu=e^{-f} dv_g$. Suppose 
that the metric and potential are time dependent, of class $\mathscr{C}^2$ and that for suitable 
constants $k, h \ge 0$ and $m\ge n$ satisfy the bounds ${\mathscr Ric}_f^m (g) \ge -(m-1)k g$ 
and $\partial_t g \ge -2 h g$ on $\mathscr M \times [t_0-T, t_0]$. Let $u$ be a positive solution 
to equation \eqref{eq11} with $1 < p < 1+1/\sqrt{m-1}$ and $v =pu^{p-1}/(p-1)$ and $M= \sup v$. 
Then there exists $C=C(p,m)>0$ such that for every $x \in \mathscr M$ and $t_0-T<t \le t_0$ we have 
\begin{align}\label{eq-2.1-EQ-global}
v^{\frac{1}{2(p-1)}}|\nabla v| (x,t) 
\le C \left \{ \begin {array}{ll}
\sqrt k M^{1+\frac{1}{2(p-1)}} + \left[ \sqrt h + \dfrac{1}{\sqrt {t-t_0+T}} \right] M^\frac{p}{2(p-1)}
\\
\\
+\sup_{\mathscr M \times [t_0-T, t_0]} \left\{\left[v^{(2p+1)/[2(p-1)]} |\Sigma_x(t,x,v)| \right]^{1/3}\right\} 
\\
\\
+ \sup_{\mathscr M \times [t_0-T, t_0]}\left\{v^{p/[2(p-1)]} \left[2 \Sigma_v(t,x,v)+\dfrac{\Sigma(t,x,v)}{(p-1)v} \right]_+^{1/2}\right\}
\end{array}
\right\}.
\end{align}
\end{theorem}

The special case of time independent metrics and potentials (the static case $\partial_t g \equiv 0$ 
and $\partial_t f \equiv 0$) is of enough significance to be formulated as a separate corollary. Here 
we describe the local version. The global version follows by passing to the limit $R\to\infty$.

\begin{theorem} \label{thm1-EQ-PME-static}
Let $(\mathscr M, g, d\mu)$ be a complete smooth metric measure space with $d\mu=e^{-f} dv_g$ and assume 
${\mathscr Ric}_f^m (g) \ge -(m-1)k g$ in ${\mathscr B}_R$ for some $k \ge 0$, $m \ge n$ 
and $R>0$. Let $u$ be a positive solution 
to \eqref{eq11} with $1 < p < 1+1/\sqrt{m-1}$ and $v =pu^{p-1}/(p-1)$ and $M= \sup_{Q_{R,T}} v$. 
Then there exists $C=C(p,m)>0$ such that for every $(x,t)$ in $Q_{R/2,T}$ with $t>t_0-T$ we have 
\begin{align}\label{eq-2.1-EQ-static}
v^{\frac{1}{2(p-1)}}|\nabla v| (x,t) 
\le C \left \{ \begin {array}{ll}
\left[ \dfrac{k^{1/4}}{\sqrt R} + \dfrac{1}{R} + \sqrt k \right] M^{1+\frac{1}{2(p-1)}}
\\
\\
+ \sup_{Q_{R, T}}\left\{v^{p/[2(p-1)]} \left[2 \Sigma_v(t,x,v)+\dfrac{\Sigma(t,x,v)}{(p-1)v} \right]_+^{1/2}\right\}
\\
\\
+ \dfrac{M^\frac{p}{2(p-1)}}{\sqrt {t-t_0+T}} 
+\sup_{Q_{R, T}} \left\{ \left[ v^{(2p+1)/[2(p-1)]} |\Sigma_x(t,x,v)| \right]^{1/3} \right\} 
\end{array}
\right\}.
\end{align}
\end{theorem}

We end this section with an application of the estimates above to parabolic Liouville-type theorems. 
Compare with Theorem \ref{ancient-1}.

\begin{theorem} \label{ancient-2}
Let $(\mathscr M, g, d\mu)$ be a complete smooth metric measure space with $d\mu=e^{-f}dv_g$ and 
${\mathscr Ric}^m_f(g) \ge 0$. Assume $(3-2p) {\mathscr N}(u) - 2u {\mathscr N}_u(u) \ge 0$ 
for all $u>0$ where $1< p < 1+1/\sqrt{m-1}$. Then any positive ancient solution to the nonlinear porous 
medium equation 
\begin{equation} \label{ancient-equation-2}
\partial_t u - \Delta_f u^p = {\mathscr N}(u(x,t)), 
\end{equation}
satisfying the growth at infinity 
\begin{equation}
u(x,t) = o \left( [\varrho(x) + \sqrt{|t|}]^{2/(2p-1)} \right), 
\end{equation}
must be spatially constant. In particular, if additionally, ${\mathscr N}(u) \ge a$ for some $a>0$ and all 
$u>0$ then \eqref{ancient-equation-2} admits no such ancient solutions. 
\end{theorem}

\section{Evolution inequalities ${\bf II}$: $\mathscr L_v^p=\partial_t - (p-1)v \Delta_f$ and $w = |\nabla v|^2/v^{\beta}$} 
\label{sec7}

In this section, similar to what was done earlier in Section \ref{sec3}, we derive evolution identities and 
inequalities for $w = |\nabla v|^2/v^{\beta}$ where the pressure $v$ relates to the positive solution $u$ 
via $v =pu^{p-1}/(p-1)$. The evolution operator here is $\mathscr L_v^p = \partial_t -(p-1) v \Delta_f$ and 
$\beta \in {\mathbb R}$ is to be specified later. Note that initially we only assume $p>1$ but later in 
the course of the proof of the estimates and bounds in the next section we further restrict the range 
of $p$. 
\begin {lemma}\label{Lem-3.1-new}
Let $(\mathscr M, g, d\mu)$ be a complete smooth metric measure space with $d\mu=e^{-f} dv_g$. 
Suppose that the metric and potential are time dependent and of class $\mathscr{C}^2$. 
Let $u$ be a positive solution to \eqref{eq11} and set 
$w = |\nabla v|^2/v^{\beta}$ where $v =pu^{p-1}/(p-1)$ and $\beta \in {\mathbb R}$ is a 
constant. Then for any constant $\varepsilon \in {\mathbb R}$, $w$ satisfies the evolution inequality 
\begin{align}\label{EQ-estimate-3.3}
\mathscr L_v^p [w]- \varepsilon \langle \nabla v, \nabla w \rangle
=& ~ [\partial_t-(p -1) v \Delta_f] w - \varepsilon \langle \nabla v, \nabla w\rangle\nonumber\\
\le &~\frac{\Gamma (\beta, \varepsilon)}{2 (p -1)} v^{\beta -1} w^2 
- 2 \left[ \frac{1}{2} \partial_t g + (p-1) v {\mathscr Ric}_f^m(g) \right] \frac{(\nabla v, \nabla v)}{v^\beta} \nonumber\\
&+2 \frac{\langle \nabla v, \Sigma_x(t,x,v) \rangle}{v^{\beta}}
+\left[2 \Sigma_v(t,x,v) -\beta \frac{ \Sigma(t,x,v)}{v}\right]w.
\end{align}
Here 
\begin{align}\label{EQ-eq-3.4-Gamma}
\Gamma (\beta, \varepsilon) =&~ [\varepsilon-2[1+(p -1) \beta] -(p -1)]^2+(m-1)(p -1)^2 \nonumber\\
&-2(p -1) \beta [1+(p -1)(\beta+1) -\varepsilon].
\end{align}
\end{lemma}

\begin{proof}
Starting from \eqref{eq2.3}  in Lemma \ref{Lem-2.3} and 
recalling $\mathscr L_v^p[w] =[\partial_t-(p -1) v \Delta_f]w$, 
$\Delta_f v  = \Delta v - \langle \nabla f, \nabla v \rangle$ we can write
\begin{align}\label{eq9.1}
\mathscr L_v^p [w] = &
-2\left[\frac{1}{2}\partial_t g+(p -1)v{\mathscr Ric}_f^m(g)\right]\frac{( \nabla v, \nabla v)}{v^\beta}
-2(p -1) \frac{|\nabla\nabla v|^2}{v^{\beta -1}}\nonumber\\
&+2(p -1) \frac{|\nabla v|^2}{v^{\beta}}  [\Delta v - \langle \nabla f, \nabla v \rangle]
-2\frac{p-1}{v^{\beta -1}}\frac{\langle \nabla f, \nabla v \rangle^2}{m-n}\nonumber\\
&+4[1+(p -1) \beta] \frac{\nabla \nabla v (\nabla v, \nabla v)}{v^\beta}
-\beta [1+(p -1)(\beta+1)]\frac{|\nabla v|^4}{v^{\beta+1}}\nonumber\\
&+2 \frac{\langle \nabla v, \Sigma_x(t,x,v) \rangle}{v^{\beta}}
+\left[2 \Sigma_v(t,x,v) -\beta \frac{ \Sigma(t,x,v)}{v}\right]\frac{|\nabla v|^2}{v^{\beta}}.
\end{align}
Now as $\langle \nabla v, \nabla w \rangle = 2 \nabla \nabla v (\nabla v, \nabla v)/v^\beta
-\beta |\nabla v|^4/v^{\beta+1}$ for any constant $\varepsilon$ we can write
\begin{align}\label{EQ-eq-3.3}
\mathscr L_v^p [w] - \varepsilon \langle \nabla v, \nabla w \rangle =
&-2(p -1) \frac{|\nabla\nabla v|^2}{v^{\beta -1}}
+2(p -1) \frac{|\nabla v|^2}{v^{\beta}}  \Delta v \\
&- \{ 2\varepsilon-4[1+(p -1) \beta] \} \frac{\nabla \nabla v (\nabla v, \nabla v)}{v^\beta}\nonumber\\
&-2\left[\frac{1}{2}\partial_t g+(p -1)v{\mathscr Ric}_f^m(g)\right]\frac{( \nabla v, \nabla v)}{v^\beta}\nonumber\\
& -\frac{2(p-1)}{v^{\beta-1}} \left[\frac{|\nabla v|^2}{v} \langle \nabla f, \nabla v \rangle
+\frac{\langle \nabla f, \nabla v \rangle^2}{(m-n)}\right] \nonumber \\
&-\beta [1+(p -1)(\beta+1) -\varepsilon ]\frac{|\nabla v|^4}{v^{\beta+1}}\nonumber\\
&+2 \frac{\langle \nabla v, \Sigma_x(t,x,v) \rangle}{v^{\beta}}
+\left[2 \Sigma_v(t,x,v) -\beta \frac{ \Sigma(t,x,v)}{v}\right]\frac{|\nabla v|^2}{v^{\beta}}.\nonumber
\end{align}

Next, focusing on the terms involving the second derivatives of $v$ on the right-hand side, 
and setting $A=\nabla \nabla v$ with ${\rm tr} A=\Delta v$ and $e = \nabla v/ |\nabla v|$, 
we can rewrite \eqref{EQ-eq-3.3} as
\begin{align}
\mbox{\eqref{EQ-eq-3.3}} 
= &- 2(p -1) \frac{|A|^2}{v^{\beta -1}} 
+2(p -1) \frac{{\rm tr} A }{|A|} w|A| 
- \{ 2\varepsilon-4[1+(p -1) \beta] \} \frac{A(e,e)}{|A|} w |A|\nonumber\\
&-2\left[\frac{1}{2}\partial_t g+(p -1) v {\mathscr Ric}_f^m(g)\right] \frac{( \nabla v, \nabla v)}{v^\beta} 
-\beta [1+(p -1)(\beta+1) -\varepsilon ]v^{\beta-1} w^2 \nonumber \\
&-\frac{2(p-1)}{v^{\beta-1}} \left[\frac{|\nabla v|^2}{v} \langle \nabla f, \nabla v \rangle
+\frac{\langle \nabla f, \nabla v \rangle^2}{(m-n)}\right] 
+2 \frac{\langle \nabla v, \Sigma_x(t,x,v) \rangle}{v^{\beta}} \nonumber \\
&+\left[2 \Sigma_v(t,x,v) -\beta \frac{ \Sigma(t,x,v)}{v}\right]w.
\end{align}
Rearranging the expression on the first line then results in 
\begin{align}\label{EQ-eq-3.5}
\mbox{\eqref{EQ-eq-3.3}} 
= &- 2(p -1) \frac{|A|^2}{v^{\beta -1}} 
-\left[ \{ 2\varepsilon-4[1+(p -1) \beta] \} \frac{A(e,e)}{|A|} 
-2(p -1) \frac{{\rm tr} A }{|A|} \right] w |A| \nonumber\\
& -2\left[\frac{1}{2}\partial_t g+(p -1) v {\mathscr Ric}_f^m(g)\right] \frac{( \nabla v, \nabla v)}{v^\beta} 
-\beta [1+(p -1)(\beta+1) -\varepsilon ] v^{\beta-1} w^2 \nonumber\\
&-\frac{2(p-1)}{v^{\beta-1}} \left[\frac{|\nabla v|^2}{v} \langle \nabla f, \nabla v \rangle
+\frac{\langle \nabla f, \nabla v \rangle^2}{(m-n)}\right]
+2 \frac{\langle \nabla v, \Sigma_x(t,x,v) \rangle}{v^{\beta}}\nonumber\\
&+\left[2 \Sigma_v(t,x,v) -\beta \frac{ \Sigma(t,x,v)}{v}\right]w.
\end{align}
Now basic algebraic manipulations involving completing squares leads to 
\begin{align}
\mathscr L_v^p& [w] - \varepsilon \langle \nabla v, \nabla w \rangle 
=[\partial_t-(p - 1) v \Delta_f]w - \varepsilon \langle \nabla v, \nabla w \rangle\nonumber\\
= &- 2(p -1) \left[ \frac{|A|}{v^{\frac{\beta-1}{2}}} + \frac{1}{2(p -1)}
\left( \{\varepsilon-2[1+(p -1) \beta] \} \frac{A(e,e)}{|A|} 
-(p -1) \frac{{\rm tr} A }{|A|} \right) v^{\frac{\beta-1}{2}} w \right]^2 \nonumber\\
& +\frac{1}{2(p -1)} \left[ \{\varepsilon-2[1+(p -1) \beta] \} \frac{A(e,e)}{|A|} 
-(p -1) \frac{{\rm tr} A }{|A|}\right]^2 v^{\beta -1} w^2\nonumber\\
& -2\left[\frac{1}{2}\partial_t g +(p -1)v{\mathscr Ric}_f^m(g)\right] \frac{( \nabla v, \nabla v)}{v^\beta} 
-\frac{2(p-1)}{v^{\beta-1}} \left[ \frac{\langle \nabla f, \nabla v \rangle}{\sqrt{m-n}} + 
\frac{|\nabla v|^2 }{2v} \sqrt{m-n}\right]^2  \nonumber\\
&+ \frac{(p -1)(m-n)}{2} v^{\beta -1} w^2
-\beta [1+(p -1)(\beta+1) -\varepsilon ] v^{\beta-1} w^2\nonumber\\
&+2 \frac{\langle \nabla v, \Sigma_x(t,x,v) \rangle}{v^{\beta}}
+\left[2 \Sigma_v(t,x,v) -\beta \frac{ \Sigma(t,x,v)}{v}\right]w.
\end{align}
Subsequently, as we are aiming for an upper bound, ignoring non-positive terms lead to the inequality
\begin{align}\label{eq9.4}
\mathscr L_v^p [w] - \varepsilon \langle \nabla v, \nabla w \rangle 
&\le \frac{1}{2 (p -1)} \left[ \{ \varepsilon-2[1+(p -1) \beta] \} \frac{A(e,e)}{|A|} 
-(p -1) \frac{{\rm tr} A}{|A|}\right]^2 v^{\beta -1} w^2\nonumber\\
&-2\left[\frac{1}{2}\partial_t g+(p -1) v {\mathscr Ric}_f^m(g)\right] \frac{( \nabla v, \nabla v)}{v^\beta} 
+ \frac{(p -1)(m-n)}{2} v^{\beta -1} w^2 \nonumber\\
&-\beta [1+(p -1)(\beta+1) -\varepsilon ] v^{\beta-1} w^2
+2 \frac{\langle \nabla v, \Sigma_x(t,x,v) \rangle}{v^{\beta}}\nonumber\\
&+\left[2 \Sigma_v(t,x,v) -\beta \frac{ \Sigma(t,x,v)}{v}\right]w.
\end{align}
Now to proceed further we make use of the following standard matrix variational result.

\begin{lemma}\label{lem9.1}
For any $a, b \in {\mathbb R}$ we have 
\begin{align}
\max_{A \in {\bf Sym}_n \atop{A \neq 0 \atop{|e|=1}}} \left[ \frac{aA + b({\rm tr} A)I_n}{|A|} (e,e) \right]^2
= (a+b)^2 +(n-1) b^2.
\end{align}
Here ${\bf Sym}_n$ is the space of real symmetric matrices and $I_n$ the identity matrix.
\end{lemma}

Applying Lemma \ref{lem9.1} to maximise the expression on the first line of inequality \eqref{eq9.4} 
(with the choices of $a= \varepsilon-2[1+(p -1) \beta]$ and $b= -(p -1)$) and rearranging terms results in  
\begin{align}\label{eq9.7}
\mathscr L_v^p [w] - \varepsilon \langle \nabla v, \nabla w \rangle 
\le &~ \frac{1}{2 (p -1)} \bigg\{ [\varepsilon-2[1+(p -1) \beta] -(p -1)]^2
+ (n-1)(p -1)^2 \nonumber\\
 &-2(p -1) \beta [1+(p -1)(\beta+1) -\varepsilon]
 + (p -1)^2(m-n)\bigg\} v^{\beta -1} w^2 \nonumber\\
&-2\left[\frac{1}{2}\partial_t g+(p -1)v{\mathscr Ric}_f^m(g)\right] \frac{( \nabla v, \nabla v)}{v^\beta}
+2 \frac{\langle \nabla v, \Sigma_x(t,x,v) \rangle}{v^{\beta}}\nonumber\\
&+\left[2 \Sigma_v(t,x,v) -\beta \frac{ \Sigma(t,x,v)}{v}\right]w.
\end{align}
Simplifying further leads to the desired conclusion.
\end{proof}

\begin {lemma}\label{Lem-3.1-new-after}
Under the assumptions of Lemma $\ref{Lem-3.1-new}$, if the metric $g$ and potential $f$ satisfy the 
super flow inequality 
\begin{align} \label{flow-inequality-SPR-PME-after}
\frac{1}{2} \partial_t g + (p-1) v {\mathscr Ric}_f^m (g) \ge - {\mathsf k} g,  
\end{align}
then $w$ satisfies the evolution inequality 
\begin{align} \label{Asasi-optimised-after}
\mathscr L_v^p [w]- p \langle \nabla v, \nabla w \rangle
\le&~ \frac{(m-1)(p -1)^2-1}{2(p -1)} \frac{w^2}{v^{p/(p-1)}} + 2{\mathsf k} w \\
&+2 v^{1/(p-1)}\langle \nabla v, \Sigma_x(t,x,v) \rangle
+\left[2 \Sigma_v(t,x,v) + \frac{ \Sigma(t,x,v)}{(p-1)v}\right]w. \nonumber 
\end{align}
\end{lemma}

Note that $\mathsf k$ here may be a constant or more generally a function of space and time. 
Moreover, for the coefficient of the term $v^{\beta-1} w^2= w^2/v^{p/(p-1)}$ (see also \eqref{eq9.10} below) 
to be negative (i.e., $<0$) it suffices to have $1<p<1+1/\sqrt{m-1}$. The significance of this 
negativity will become apparent later ({\it see}, e.g., \eqref{EQ-eq-2.20} in the next section).

\begin{proof}
Optimising the quadratic function $\Gamma (\beta, \varepsilon)$ defined by \eqref{EQ-eq-3.4-Gamma} 
in the inequality \eqref{EQ-estimate-3.3} leads to the optimiser $(\beta^\star, \varepsilon^\star)$ with 
\begin{align}\label{pre-eq9.10}
\beta^\star = -\frac{1}{p-1}, \qquad \text{and} \qquad \varepsilon^\star = p.
\end{align}
Substituting these into \eqref{EQ-eq-3.4-Gamma} then leads to the optimal value 
\begin{align}\label{eq9.10}
\Gamma(\beta^\star, \varepsilon^\star) = \Gamma (-1/(p-1), p) = (m-1)(p-1)^2-1.
\end{align}
Using the optimal values \eqref{pre-eq9.10} and \eqref{eq9.10} in the inequality \eqref{EQ-estimate-3.3} 
with $w=|\nabla v|^2 /v^\beta = v^{1/(p-1)}|\nabla v|^2$ then results in
%and using the bound $(1/2) \partial_t g + (p-1) v {\mathscr Ric}_g^m(g) \ge -[(p-1)(m-1) kv+h]g$ 
%(in the sense of symmetric tensors) results in the inequality
\begin{align} \label{Asasi-optimised}
\mathscr L_v^p [w]- p \langle \nabla v, \nabla w \rangle
\le&~ \frac{(m-1)(p -1)^2-1}{2(p -1)} \frac{w^2}{v^{p/(p-1)}} \nonumber \\
& - 2 v^{1/(p-1)} \left[ \frac{1}{2} \partial_t g+(p -1) v {\mathscr Ric}_f^m(g) \right] ( \nabla v, \nabla v) \nonumber\\
&+2 v^{1/(p-1)}\langle \nabla v, \Sigma_x(t,x,v) \rangle
+\left[2 \Sigma_v(t,x,v) + \frac{ \Sigma(t,x,v)}{(p-1)v}\right]w.\nonumber
\end{align}
Invoking the super flow inequality \eqref{flow-inequality-SPR-PME-after} gives at once the desired conclusion. 
\end{proof}

\section{Cylindrical localisation and the proof of Theorem \ref{thm1-EQ-PME}}
\label{sec8}

Starting again with the localised function $\eta w$, where $\eta$ is as in \eqref{cut-off def}, it follows upon using the identity in Lemma \ref{Lem-2.6} that  
\begin{align}
\mathscr L_v^p[\eta w] = w\mathscr L_v^p[\eta]
- 2(p -1)v [\langle\nabla \eta ,\nabla (\eta w) \rangle- |\nabla \eta|^2 w]/\eta
+\eta \mathscr L_v^p [w].
\end{align}
Using \eqref{Asasi-optimised} along with the observation that 
$\eta \langle \nabla v, \nabla w \rangle = \langle \nabla v, \nabla (\eta w) \rangle - w \langle \nabla v, \nabla \eta \rangle$ then gives
\begin{align}\label{EQ-eq-8.18}
\mathscr L_v^p [\eta w] \le &~
w \mathscr L_v^p [\eta]
-2(p -1) v \left\langle \frac{\nabla \eta}{\eta} , \nabla (\eta w) \right\rangle
+ 2(p -1)v \frac{|\nabla \eta|^2}{\eta} w \nonumber\\
&+\varepsilon^\star \langle \nabla v , \nabla (\eta w) \rangle
-\varepsilon^\star w\langle \nabla v , \nabla \eta \rangle
+\frac{\Gamma (\beta^\star, \varepsilon^\star)}{2 (p -1)} v^{\beta^\star -1} \eta w^2 \nonumber\\
&+ 2[(p - 1)(m-1)k v+h] \eta w
+ 2\eta \frac{\langle \nabla v, \Sigma_x(t,x,v) \rangle}{v^{\beta^\star}} \nonumber\\
&+ \left[2 \Sigma_v(t,x,v) -\frac{ \beta^\star \Sigma(t,x,v)}{v} \right] \eta w.
\end{align}  
Assume that the localised function $\eta w$ is maximised at the point $(x_1, t_1)$ in the compact space-time set 
$\{d(x,x_0, t) \le R, t_0-T \le t \le \tau\}$. Additionally, we can assume that $(\eta w)(x_1, t_1) >0$ 
as otherwise the conclusion of the theorem is true with $w(x, \tau) \le 0$ for all $d(x, x_0, \tau) \le R/2$. Hence 
$t_1>t_0-T$ and at the point $(x_1,t_1)$  we have $\Delta_f(\eta w) \le 0$, $\partial_t (\eta w) \ge 0$ and 
$\nabla(\eta w) =0$. Thus at the point $(x_1,t_1)$ we also have $\mathscr L_v^p [\eta w] \ge 0$. 
From \eqref{EQ-eq-8.18} upon rearranging of terms it therefore follows that  
\begin{align}
- \frac{\Gamma (\beta^\star, \varepsilon^\star)}{2 (p -1)} v^{\beta^\star-1} \eta w^2 \le
&~ w \mathscr L_v^p [\eta]
+2[(p - 1)(m-1) k v + h] \eta w\nonumber\\
&- \varepsilon^\star w\langle \nabla v , \nabla \eta \rangle 
+ 2(p -1)v \frac{|\nabla \eta|^2}{\eta} w\nonumber\\
&+\left[2 \Sigma_v(t,x,v)-\frac{ \beta^\star \Sigma(t,x,v)}{v} \right] \eta w %\nonumber\\
%&
+ 2\eta \frac{\langle \nabla v, \Sigma_x(t,x,v) \rangle}{v^{\beta^\star}}.
\end{align}

Denoting $-\Gamma (\beta^\star, \varepsilon^\star)/[2(p -1)]= 2/\gamma>0$, that is, $\gamma=4(p-1)/[1-(m-1) (p-1)^2]$ 
(see the comments following Lemma \ref {Lem-3.1-new-after}) and multiplying through by $\gamma v^{1-\beta^\star}$ gives 
\begin{align}\label{EQ-eq-2.20}
2\eta w^2 \le
&~ \gamma v^{1-\beta^\star} w \mathscr L_v^p [\eta]
+2 \gamma [(p - 1)(m-1) k_1 v + k_2]v^{1-\beta^\star} \eta w\nonumber\\
&- \varepsilon^\star \gamma v^{1-\beta^\star} w\langle \nabla v , \nabla \eta \rangle 
+ 2 \gamma (p -1)v^{2-\beta^\star} \frac{|\nabla \eta|^2}{\eta} w\nonumber\\
&+\gamma v^{1-\beta^\star} \eta w \left[2 \Sigma_v(t,x,v)-\frac{ \beta^\star \Sigma(t,x,v)}{v} \right] %\nonumber\\
%&
+ 2 \gamma \eta v^{1-\beta^\star} \frac{\langle \nabla v, \Sigma_x(t,x,v) \rangle}{v^{\beta^\star}}.
\end{align}

Now we proceed onto bounding the expression on the right-hand side of \eqref{EQ-eq-2.20}.
As this is similar to the proof of Theorem \ref{thm1} we shall remain brief, focusing mainly on the differences.

Towards this end, starting from the first term on the right and arguing as in Theorem \ref{thm1} we have 
\begin{align}
\gamma v^{1-\beta^\star} w \mathscr L_v^p [\eta] \le&~ \frac{1}{7} \eta w^2
+ C \gamma ^2 \left[\frac{1}{(\tau-t_0+T)^2} + h^2\right](\sup_{Q_{R, T}} v)^{2-2\beta^\star}\nonumber\\
&+ C \gamma^2(p -1)^2 (m-1)^2\left[\frac{1+ k R^2}{R^4}\right] (\sup_{Q_{R, T}} v)^{4-2\beta^\star}.
\end{align}
Likewise for the second term we can write 
\begin{align}
2 \gamma (p - 1)(m-1) k v^{2-\beta^\star} \eta w 
\le \frac{1}{7} \eta w^2 + C \gamma ^2 (p-1)^2 (m-1)^2 k^2(\sup_{Q_{R, T}} v)^{4-2\beta^\star}, \\
2 \gamma h v^{1-\beta^\star} \eta w \le
2 \gamma h v^{1-\beta^\star} \eta^{1/2} w 
\le \frac{1}{7} \eta w^2 + C \gamma ^2 h^2 (\sup_{Q_{R, T}} v)^{2-2\beta^\star}.
\end{align} 
We can bound the next term by writing 
\begin{align}\label{EQ-bound1}
-\varepsilon^\star \gamma v^{1-\beta^\star}w \langle \nabla v , \nabla \eta \rangle 
& \le \varepsilon^\star \gamma v^{1-\beta^\star}w |\nabla v| |\nabla \eta| \nonumber\\
& \le \varepsilon^\star \gamma v^{1-\beta^\star/2} \eta^{3/4} w^{3/2} \frac{|\nabla \eta|}{\eta^{3/4}}\nonumber\\
& \le \frac{1}{7} (\eta ^{3/4} w^{3/2})^{4/3} + C \varepsilon^{\star 4} \gamma^4 \left[\frac{|\nabla \eta|}{\eta^{3/4}}\right]^4 v^{4-2\beta^\star} \nonumber\\
& \le \frac{1}{7} \eta w^2 + \frac{C \varepsilon^{\star 4} \gamma^4}{R^4} (\sup_{Q_{R, T}} v)^{4-2\beta^\star}.
\end{align}
In much the same way
\begin{align}
2 \gamma (p -1) v^{2-\beta^\star} \frac{|\nabla \eta|^2}{\eta} w 
\le \frac{1}{7} \eta w^2 + \frac{C \gamma^2 (p -1)^2}{R^4}(\sup_{Q_{R, T}} v)^{4-2\beta^\star}.
\end{align}
 
Finally for the $\Sigma(t,x,v)$ terms, upon recalling $0\le \eta \le 1$ and utilising Young's inequality we have
\begin{align}\label{eq-2.33}
\gamma v^{1-\beta^\star}\eta w \bigg[2 \Sigma_v(t,x,v)&-\frac{ \beta^\star \Sigma(t,x,v)}{v} \bigg]\\
& \le \frac{1}{7} \eta w^2 + C \gamma^2 \sup_{Q_{R, T}} 
\left\{v^{2(1-\beta^\star)}\left[2 \Sigma_v(t,x,v)-\frac{ \beta^\star \Sigma(t,x,v)}{v} \right]_+^2\right\}.\nonumber
\end{align}
Furthermore, by using Cauchy-Schwarz and Young inequalities we can write
\begin{align}\label{eq-2.32}
2 \gamma v^{1-\beta^\star} \eta \frac{\langle \nabla v, \Sigma_x(t,x,v) \rangle}{v^{\beta^\star}} \le &~
2 \gamma v^{1-\beta^\star} \eta |\nabla v| \frac{|\Sigma_x(t,x,v)|}{v^{\beta^\star}}\\
\le &~ \frac{1}{7} \eta w^2+C \gamma^{4/3} \sup_{Q_{R, T}} 
\left\{\left[\frac{|\Sigma_x(t,x,v)|}{v^{(3\beta^\star-2)/2}}\right]^{4/3}\right\}.\nonumber
\end{align}

Completing the estimate of the individual terms on the right-hand side of \eqref{EQ-eq-2.20} 
we now proceed by substituting these back into the inequality. Recalling the explicit optimal 
values \eqref{pre-eq9.10} and using the local bound $M=\sup_{Q_{R,T}} v$ it then follows from 
\eqref{EQ-eq-2.20} and the individual bounds above (note that the coefficient $2$ of $\eta w^2$ 
on the left-hand side exceeds the sum $1$ of the coefficients of the seven term involving 
$\eta w^2$ on the right-hand side where each coefficient appears as $1/7$) that
\begin{align}
\eta w^2 \le C(\gamma) &
\left\{ \left( \begin {array}{ll} 
(p -1)^2 (m-1)^2\left[\dfrac{1+ k R^2}{R^4}\right] +\dfrac{p^4}{R^4} 
\\
\\
+\dfrac{(p -1)^2}{R^4} + (p-1)^2 (m-1)^2 k^2
\end{array}
\right) \right. M^{4+2/(p-1)}
\nonumber \\
&\qquad + 
\left. \left(  \begin{array}{ll}
\left[ 2h^2 + \dfrac{1}{(\tau-t_0+T)^2} \right] M^{2+2/(p-1)}
\\
\\
\sup_{Q_{R, T}} \left\{v^{2p/(p-1)} \left[2 \Sigma_v(t,x,v)+\dfrac{\Sigma(t,x,v)}{(p-1)v} \right]_+^2\right\}
\\
\\
+ \sup_{Q_{R, T}} \left\{[v^{(2p+1)/[2(p-1)]}|\Sigma_x(t,x,v)|]^{4/3}\right\}
\end{array}
\right) \right\}.
\end{align}

Recalling the maximality of $\eta w$ at $(x_1, t_1)$ along with $\eta \equiv 1$ when $d(x, x_0,t) \le R/2$ and $\tau \le t \le t_0$, 
it follows that $w^2 (x, \tau) = (\eta^2 w^2)(x, \tau) \le (\eta^2 w^2)(x_1,t_1) \leq (\eta w^2)(x_1,t_1)$. 
Hence noting $w = v^{1/(p-1)}|\nabla v|^2$ and 
absorbing $p$ and $m$ into $C(\gamma)$ gives
\begin{align}
[v^{1/[2(p-1)]}|\nabla v|] (x,\tau)
\le C \left \{ \begin {array}{ll}
\sqrt h M^\frac{p}{2(p-1)}
+ \left[ \dfrac{k^{1/4}}{\sqrt {R}} 
+\dfrac{1}{R} +\sqrt k \right] M^{1+\frac{1}{2(p-1)}}
\\
\\
+ \sup_{Q_{R, T}}\left\{v^\frac{p}{2(p-1)} \left[2 \Sigma_v(t,x,v)+\dfrac{\Sigma(t,x,v)}{(p-1)v} \right]_+^{1/2}\right\}
\\
\\
+ \dfrac{M^\frac{p}{2(p-1)}}{\sqrt {\tau-t_0+T}} 
+ \sup_{Q_{R, T}} \left\{[v^\frac{2p+1}{2(p-1)} |\Sigma_x(t,x,v)|]^{1/3}\right\}
\end{array}
\right\}.
\end{align}
The arbitrariness of $t_0-T<\tau \le t_0$ now gives the desired conclusion.
\hfill $\square$

\section{Ancient solutions to $\partial_t u - \Delta_f u^p = {\mathscr N}(u)$ ${\bf II}$: Proof of Theorem \ref{ancient-2}}
\label{sec9}

Let us now present the proof of Theorem \ref{ancient-2}. We first fix a space-time point $(x_0,t_0)$. 
Then with the choices $t=t_0$, $R>0$ and $T=R^2$ it follows from the growth assumption 
$u(x,t) = o([\varrho(x) + \sqrt{|t|}]^{2/(2p-1)})$ that for $R$ sufficiently large 
\begin{equation} \label{M-growth-ancient-2}
M = \sup_{Q_{R,T}} v = [p/(p-1)] \sup_{Q_{R,T}} u^{p-1} = o(R^\frac{p-1}{p-1/2}).
\end{equation} 
Now turning to the local estimate \eqref{eq-2.1-EQ-static} in Theorem \ref{thm1-EQ-PME-static} 
(with $k=0$) we can write 
\begin{align} \label{vx0t0-ancient-2}
v^{\frac{1}{2(p-1)}}|\nabla v| (x_0,t_0) 
&\le C \left \{ \begin {array}{ll}
\left[ \dfrac{k^{1/4}}{\sqrt {R}} 
+\dfrac{1}{R} +\sqrt k \right] M^{1+\frac{1}{2(p-1)}}
\\
\\
+ \sup_{Q_{R, T}}\left\{v^\frac{p}{2(p-1)} \left[2 \Sigma_v(t,x,v)+\dfrac{\Sigma(t,x,v)}{(p-1)v} \right]_+^\frac{1}{2}\right\}
\\
\\
+ \dfrac{M^\frac{p}{2(p-1)}}{\sqrt {\tau-t_0+T}} 
+ \sup_{Q_{R, T}} \left\{[v^\frac{2p+1}{2(p-1)} |\Sigma_x(t,x,v)|]^\frac{1}{3}\right\}
\end{array}
\right\} \\
&\le C \left \{ \begin {array}{ll}
\dfrac{M^{1+\frac{1}{2(p-1)}}}{R} 
+\sup_{Q_{R, T}} \left\{ v^\frac{2p+1}{2(p-1)} \left[ |\Sigma_x(t,x,v)| \right]^\frac{1}{3} \right\} 
\\
\\
\dfrac{M^{\frac{p}{2(p-1)}} }{\sqrt T} + \sup_{Q_{R, T}}\left\{v^\frac{p}{2(p-1)} \left[2 \Sigma_v(t,x,v)+\dfrac{\Sigma(t,x,v)}{(p-1)v} \right]_+^\frac{1}{2}\right\}
\end{array}
\right\}. \nonumber 
\end{align}

Next referring to Remark \ref{Sigma-N-u-v-remark} we have $\Sigma(v) = p u^{p-2} {\mathscr N}(u)$ 
where $u=[(p-1)v/p]^{1/(p-1)}$. Thus in particular $\Sigma_x \equiv 0$ and $\Sigma_v=\Sigma_u \partial_v u =(p-2){\mathscr N}/u+{\mathscr N}_u$. 
Hence 
\begin{align*}
2 \Sigma_v(v) + \frac{\Sigma(v)}{(p-1)v} 
&= 2 \left[ (p-2) \frac{{\mathscr N}(u)}{u} + {\mathscr N}_u(u) \right] + \frac{{\mathscr N}(u)}{u} \nonumber \\
&= [2(p-2) +1] \frac{{\mathscr N}(u)}{u} + 2 {\mathscr N}_u(u) \le 0,
\end{align*}
where the last inequality follows from the assumptions on ${\mathscr N}$ in the theorem. 
Hence by virtue of \eqref{M-growth-ancient-2} we deduce from \eqref{vx0t0-ancient-2} that 
\begin{align}
v^{\frac{1}{2(p-1)}}|\nabla v| (x_0,t_0) \le \dfrac{M^{1+\frac{1}{2(p-1)}} }{R} + \dfrac{M^{\frac{p}{2(p-1)}} }{\sqrt T} 
\le \frac{o(R)}{R} + \frac{o(R^\frac{p}{2p-1})}{R}.
\end{align}
Passing to the limit $R \nearrow \infty$ and noting that $p/(2p-1) \le 1$ in the range $p>1$ 
it follows that $|\nabla v|(x_0,t_0)=0$. The arbitrariness of $(x_0,t_0)$ implies $|\nabla v| \equiv 0$ 
and so $v$ and subsequently $u$ are spatially constant. Hence we have $u=u(t)$ and from equation 
\eqref{ancient-equation-2} it follows that $du/dt = {\mathscr N}(u)$. The rest is similar to the proof of Theorem \ref{ancient-1}. 
 \hfill $\square$

\section{A Ricci-Perelman super flow inequality $\partial_t g + 2(p-1) v {\mathscr Ric}^m_f(g) \ge - 2\mathsf{k}g$} 
\label{sec10}

In this final section of the paper we derive a global gradient bound on positive smooth solutions 
to equation \eqref{eq11}. Here $\mathscr M$ is closed and the metric and potential are assumed to evolve 
under the super flow inequality 
\begin{align} \label{eq11c-2.n} 
\begin{cases}
\dfrac{1}{2} \dfrac{\partial g}{\partial t} (x,t) + (p-1) v(x,t) {\mathscr Ric}^m_f(g)(x,t) 
\ge - \mathsf{k}g(x,t), 
\\
\\ 
{\mathscr Ric}_f^m(g)(x,t) = {\mathscr Ric} (g)(x,t) + \nabla_{g} \nabla_{g} f (x,t) \\
\qquad \qquad \qquad \quad - \dfrac{\nabla_{g} f(x,t) \otimes \nabla_{g} f(x,t)}{m-n}.    
\end{cases}
\end{align}
Here $p>1$ but we will refine this range further later. (See Corollaries \ref{first-cor-last} 
and \ref{second-cor-last}.)

\begin{lemma} \label{evolution Rpq-2.n}
Let $u$ be a positive smooth solution to \eqref{eq11} and set $v =pu^{p-1}/(p-1)$. For $s \ge 2$, 
$q \in {\mathbb R}$ and $\zeta=\zeta(t)$ non-negative and of class ${\mathscr C}^1[0, \infty)$ let 
\begin{equation} \label{Xpq def-2.n}
{\mathsf H}^{s,q}_\zeta[v] = \zeta(t) \frac{|\nabla v|^s}{v^q} + \Gamma(v),
\end{equation} 
where $\Gamma=\Gamma(v)$ with $v>0$ is of class $\mathscr{C}^2$.
Then, for every $\varepsilon \in {\mathbb R}$, ${\mathsf H} = {\mathsf H}^{s, q}_\zeta [v]$ satisfies 
the evolution identity,  
\begin{align} \label{R pq phi equality-2}
{\mathscr L}_v^{p} ({\mathsf H}^{s, q}_\zeta [v]) 
&- \varepsilon \langle \nabla v, \nabla {\mathsf H}^{s, q}_\zeta [v] \rangle
= [\partial _t -(p-1) v \Delta_ f] ({\mathsf H}^{s, q}_\zeta [v]) 
- \varepsilon \langle \nabla v, \nabla {\mathsf H}^{s, q}_\zeta [v] \rangle\nonumber\\
=&~ \zeta'(t) \frac{|\nabla v|^s}{v^q} - s \zeta (t) \frac{|\nabla v|^{s-2}}{v^q} 
\left[ \frac{1}{2} \partial_t g + (p-1) v {\mathscr Ric}_f^m (g)\right] (\nabla v, \nabla v) \nonumber \\
&+ \zeta (t) \frac{|\nabla v|^{s-2}}{v^q} \left\{ s (p-1)  \left[|\nabla v|^2 \Delta_f v - v |\nabla \nabla v|^2
- v \frac{\langle \nabla f , \nabla v \rangle^2}{(m-n)} \right] \right. \nonumber\\
&\left. + \left[ s[q(p-1)+1-\varepsilon/2] \langle \nabla |\nabla v|^2, \nabla v\rangle
- q[q(p-1) +p-\varepsilon] \frac{|\nabla v|^{4}}{v}\right] \right\} \nonumber\\
&- \zeta (t)\frac{s(p-1)}{2v^{q-1}} \langle \nabla |\nabla v|^{s-2} , \nabla |\nabla v|^2 \rangle
+ s \zeta (t) \frac{|\nabla v|^{s-2}}{v^q}\langle \nabla v, \nabla \Sigma(t,x,v) \rangle\nonumber\\
& - q \zeta (t)\frac{|\nabla v|^s}{v^q} \frac{\Sigma(t,x,v)}{v} 
- (p-1)v\Gamma''(v)  |\nabla v|^2\nonumber\\
&+ \Gamma'(v) [(1-\varepsilon)|\nabla v|^2 +\Sigma(t,x,v)].
\end{align}
\end{lemma}

\begin{proof}
Referring to \eqref{Xpq def-2.n} it is a straightforward matter to see that  
\begin{align}\label{eq-10.5-2.n}
{\mathscr L}_v^{p} ({\mathsf H}^{s,q}_\zeta[v]) = \zeta'(t) \frac{|\nabla v|^s}{v^q} 
+ \zeta (t) {\mathscr L}_v^{p} \left[\frac{|\nabla v|^s}{v^q} \right] + {\mathscr L}_v^{p} [\Gamma(v)].
\end{align}

We need to evaluate the expressions in the second and third terms on the right-hand side respectively. 
Focusing on the second term, first it is easily seen that,  
\begin{align}\label{10.6.n-2.n}
{\mathscr L}_v^{p} \left[\frac{|\nabla v|^s}{v^q} \right] 
=&~ \frac{1}{v^q} {\mathscr L}_v^{p} [|\nabla v|^s] 
+ 2 (p-1) v \left \langle \nabla \left(\frac{|\nabla v|^s}{v^q}\right), \frac{\nabla v^q}{v^q} \right \rangle 
- \frac{|\nabla v|^s}{v^{2q}}  {\mathscr L}_v^{p} [v^q].
\end{align}
Now for the first term on the right-hand side of \eqref{10.6.n-2.n} a direct calculation gives  
\begin{align}\label{EQ-eq10.7.n-2.n}
{\mathscr L}_v^{p} [|\nabla v|^s] =&~  [\partial _t -(p-1) v \Delta_ f] |\nabla v|^s 
= \partial _t (|\nabla v|^2)^{s/2} -(p-1) v \Delta_ f (|\nabla v|^2)^{s/2} \nonumber\\
%=&~ \frac{s}{2} |\nabla v|^{s-2} \partial_t |\nabla v|^2 - \frac{s}{2}(p-1) v |\nabla v|^{s-2} \Delta_f |\nabla v|^2 
%- \frac{s}{2}(p-1)v \langle \nabla |\nabla v|^{s-2} , \nabla |\nabla v|^2 \rangle \nonumber\\
=&~\frac{s}{2} |\nabla v|^{s-2} \left[ \partial_t -(p-1)v \Delta_f\right] |\nabla v|^2
- \frac{s}{2}(p-1)v \langle \nabla |\nabla v|^{s-2} , \nabla |\nabla v|^2 \rangle \nonumber\\
=&~ \frac{s}{2} |\nabla v|^{s-2} {\mathscr L}_v^{p} [|\nabla v|^2] 
- \frac{s}{2}(p-1)v \langle \nabla |\nabla v|^{s-2} , \nabla |\nabla v|^2 \rangle, 
\end{align}
whilst using Lemma \ref{Lem.2.2} and the weighted Bochner-Weitzenb\"ock formula \eqref{Bochner-1} we have  
\begin{align}
\mathscr L^p_v [|\nabla v|^2]
= & -[\partial_t g] ( \nabla v, \nabla v) +2 (p-1) |\nabla v|^2 \Delta_f v + 2\langle \nabla v, \nabla |\nabla v|^2 \rangle\nonumber\\
&+2 \langle \nabla v, \nabla \Sigma(t,x,v) \rangle-2(p-1) v |\nabla \nabla v|^2\nonumber\\
& -2(p-1)v {\mathscr Ric}_f^m (\nabla v, \nabla v) - 2 (p-1) v \langle \nabla f , \nabla v \rangle^2/(m-n). 
\end{align}
Thus combining the two identities above it follows that   
\begin{align}\label{10.9.n-2.n}
{\mathscr L}_v^{p} [|\nabla v|^s] =&~ \frac{s}{2} |\nabla v|^{s-2} {\mathscr L}_v^{p} [|\nabla v|^2] 
- \frac{s}{2}(p-1)v \langle \nabla |\nabla v|^{s-2} , \nabla |\nabla v|^2 \rangle\nonumber\\
=&-\frac{s}{2} |\nabla v|^{s-2} [\partial_t g] ( \nabla v, \nabla v) +  s (p-1) |\nabla v|^s \Delta_f v
+ s  |\nabla v|^{s-2}\langle \nabla v, \nabla |\nabla v|^2 \rangle\nonumber\\
&+ s |\nabla v|^{s-2}\langle \nabla v, \nabla \Sigma(t,x,v) \rangle- s (p-1) v |\nabla v|^{s-2} |\nabla \nabla v|^2\nonumber\\
& - s (p-1)v |\nabla v|^{s-2} {\mathscr Ric}_f^m (\nabla v, \nabla v)
-s (p-1)v |\nabla v|^{s-2} \langle \nabla f , \nabla v \rangle^2/(m-n)\nonumber\\
&-\frac{s}{2}(p-1)v \langle \nabla |\nabla v|^{s-2} , \nabla |\nabla v|^2 \rangle.
\end{align}
Next, for the second term on the right-hand side of \eqref{10.6.n-2.n} we can write 
\begin{align}\label{10.10.n-2.n}
\left \langle \nabla \left(\frac{|\nabla v|^s}{v^q}\right), \frac{\nabla v^q}{v^q} \right \rangle
%=&~ \frac{q}{v^{q+1}} \langle \nabla |\nabla v|^s , \nabla v\rangle  - q^2 \frac{|\nabla v|^{s+2}}{v^{q+2}}\nonumber\\
=&~\frac{sq}{2} \frac{|\nabla v|^{s-2}}{v^{q+1}}\langle \nabla |\nabla v|^2 , \nabla v\rangle  - q^2 \frac{|\nabla v|^{s+2}}{v^{q+2}}. 
\end{align}
In a similar way for the last term on the right-hand side of \eqref{10.6.n-2.n} we have
\begin{align}\label{EQ-eq-1011.n-2.n}
{\mathscr L}_v^{p} [v^q] 
&=q v^{q-1} \partial_t v -(p-1) v q v^{q-1} \Delta_f v 
-(p-1) q(q-1) v^{q-1} |\nabla v|^2\nonumber\\
&=q v^{q-1} [\partial_t - (p-1) v \Delta_f] v -(p-1) q(q-1) v^{q-1} |\nabla v|^2\nonumber\\
&=q v^{q-1} {\mathscr L}_v^{p} [v] -(p-1) q(q-1) v^{q-1} |\nabla v|^2 \nonumber \\
&=q v^{q-1} |\nabla v|^2+ q v^{q-1}\Sigma(t,x,v) -(p-1)q(q-1) v^{q-1} |\nabla v|^2,
\end{align}
where in concluding the last line we have used \eqref{eq-2.9}. 
Substituting \eqref{10.9.n-2.n}, \eqref{10.10.n-2.n} and \eqref{EQ-eq-1011.n-2.n} back into \eqref{10.6.n-2.n} 
and rearranging terms now leads to 
\begin{align}\label{eq-10.14-2.n}
{\mathscr L}_v^{p} \left[\frac{|\nabla v|^s}{v^q} \right] 
=& - s \frac{|\nabla v|^{s-2}}{v^q} \left[ \frac{1}{2} \partial_t g + (p-1) v {\mathscr Ric}_f^m(g) \right] (\nabla v, \nabla v) 
 - \frac{q}{v^{q+1}} |\nabla v|^s\Sigma(t,x,v) \nonumber\\
&+ s (p-1) \frac{|\nabla v|^{s-2}}{v^q} \left[|\nabla v|^2 \Delta_f v - v |\nabla \nabla v|^2
- v \frac{\langle \nabla f , \nabla v \rangle^2}{(m-n)} \right] \\
&+ s \frac{|\nabla v|^{s-2}}{v^{q-1}} \left[\frac{1}{v} [q(p-1)+1] \langle \nabla |\nabla v|^2 , \nabla v\rangle
- \frac{q}{s}[q(p-1) +p] \frac{|\nabla v|^{4}}{v^2}\right]\nonumber\\
&-\frac{s}{2v^{q-1}}(p-1) \langle \nabla |\nabla v|^{s-2} , \nabla |\nabla v|^2 \rangle
+ s \frac{|\nabla v|^{s-2}}{v^q}\langle \nabla v, \nabla \Sigma(t,x,v) \rangle. \nonumber
\end{align}

This completes the calculation of the second term on the right-hand side of  \eqref{eq-10.5-2.n}. Since for the last 
term of the same equation we can write 
\begin{align}\label{eq-10.15-2.n}
{\mathscr L}_v^{p} [\Gamma(v)] %&= \Gamma'(v) {\mathscr L}_v^p(v) - (p-1)v\Gamma''(v)  |\nabla v|^2\nonumber\\
&= \Gamma'(v) [|\nabla v|^2 + \Sigma(t,x,v)] - (p-1)v\Gamma''(v)  |\nabla v|^2, 
\end{align}
it follows at once upon substituting \eqref{eq-10.14-2.n} and \eqref{eq-10.15-2.n} back into \eqref{eq-10.5-2.n} that  
\begin{align}\label{eq-10.16-2}
{\mathscr L}_v^{p} ({\mathsf H}^{s,q}_\zeta[v]) 
=&~ \zeta'(t) \frac{|\nabla v|^s}{v^q} - s \zeta (t) \frac{|\nabla v|^{s-2}}{v^q} 
\left[ \frac{1}{2} \partial_t g + (p-1) v {\mathscr Ric}_f^m(g) \right] (\nabla v, \nabla v) \nonumber \\
&+\zeta (t) s (p-1)\frac{|\nabla v|^{s-2}}{v^q} \left[|\nabla v|^2 \Delta_f v - v |\nabla \nabla v|^2
- v \frac{\langle \nabla f , \nabla v \rangle^2}{(m-n)} \right] \\
&+ \zeta (t) \frac{|\nabla v|^{s-2}}{v^{q}} \left[ s[q(p-1)+1] \langle \nabla |\nabla v|^2 , \nabla v\rangle
- q[q(p-1) +p] \frac{|\nabla v|^{4}}{v}\right]\nonumber\\
&- \zeta (t)\frac{s}{2v^{q-1}}(p-1) \langle \nabla |\nabla v|^{s-2} , \nabla |\nabla v|^2 \rangle
+ \zeta (t) s \frac{|\nabla v|^{s-2}}{v^q}\langle \nabla v, \nabla \Sigma(t,x,v) \rangle\nonumber\\
& - \zeta (t)\frac{q}{v^{q+1}} |\nabla v|^s\Sigma(t,x,v)+ \Gamma'(v) [|\nabla v|^2 +\Sigma(t,x,v)]
- (p-1)v\Gamma''(v)  |\nabla v|^2. \nonumber 
\end{align}
Finally since   
\begin{align}
\left \langle \nabla v, \nabla \left(\frac{|\nabla v|^s}{v^q}\right) \right\rangle 
= \frac{s}{2} \frac{|\nabla v|^{s-2}}{v^q} \langle \nabla v, \nabla |\nabla v|^2 \rangle
- q \frac{|\nabla v|^{s+2}}{v^{q+1}},
\end{align} 
using \eqref{eq-10.16-2} and the above it follows that for any constant $\varepsilon$ we have
\begin{align}\label{eq-15.17.n-2}
{\mathscr L}_v^{p}({\mathsf H}^{s,q}_\zeta[v]) 
&- \varepsilon \zeta(t)\langle \nabla v, \nabla (|\nabla v|^s/v^q) \rangle \nonumber \\
=&~ \zeta'(t) \frac{|\nabla v|^s}{v^q} - s \zeta (t) \frac{|\nabla v|^{s-2}}{v^q} 
\left[ \frac{1}{2} \partial_t g + (p-1) v {\mathscr Ric}_f^m(g) \right] ( \nabla v, \nabla v) \\
&+\zeta (t) \frac{|\nabla v|^{s-2}}{v^q} \left\{ s (p-1) \left[|\nabla v|^2 \Delta_f v - v |\nabla \nabla v|^2
- v \frac{\langle \nabla f , \nabla v \rangle^2}{(m-n)} \right] \right. \nonumber\\
& \left. + \bigg[ s[q(p-1)+1-\varepsilon/2] \langle \nabla |\nabla v|^2 , \nabla v\rangle
- q[q(p-1) +p-\varepsilon] \frac{|\nabla v|^{4}}{v}\bigg] \right\} \nonumber\\
&- \zeta (t)\frac{s}{2v^{q-1}}(p-1) \langle \nabla |\nabla v|^{s-2} , \nabla |\nabla v|^2 \rangle
+ \zeta (t) s \frac{|\nabla v|^{s-2}}{v^q}\langle \nabla v, \nabla \Sigma(t,x,v) \rangle\nonumber\\
& - \zeta (t)\frac{q}{v^{q+1}} |\nabla v|^s\Sigma(t,x,v)+ \Gamma'(v) [|\nabla v|^2 +\Sigma(t,x,v)]
- (p-1)v\Gamma''(v)  |\nabla v|^2. \nonumber 
\end{align}
This after a rearrangement of terms gives the desired conclusion. 
\end{proof}

\begin{lemma} \label{evolution Rpq-2.n-inequality}
Under the assumptions of Lemma $\ref{evolution Rpq-2.n}$, if the metric $g$ and potential $f$ evolve under 
the flow inequality \eqref{eq11c-2.n} then 
\begin{align} \label{EQ-eq-4.4-2.n}
{\mathscr L}_v^{p} ({\mathsf H}^{s,q}_\zeta[v]) 
- &~\varepsilon \langle \nabla v, \nabla {\mathsf H}^{s,q}_\zeta[v] \rangle \nonumber\\
\le& \left\{ \zeta'(t) + s \mathsf{k} \zeta(t) + \zeta(t) \left[s\Sigma_v(t,x,v) 
- q \frac{\Sigma(t,x,v)}{v} \right] \right\} \frac{|\nabla v|^s}{v^q} \nonumber\\
&+ \frac{s\zeta (t)}{4(p-1)} \bigg[\{\varepsilon- 2[q(p-1)+1] -(p-1)\}^2 + (p-1)^2(m-1) \nonumber\\
&- \frac{4(p-1)}{s} q[q(p-1) +p-\varepsilon]\bigg] \frac{|\nabla v|^{s+2}}{v^{q+1}}
+ s \zeta (t) \frac{|\nabla v|^{s-2}}{v^q}\langle \nabla v, \Sigma_x (t,x,v)\rangle  \nonumber\\
&- (p-1)v\Gamma''(v)  |\nabla v|^2 + \Gamma'(v) [(1-\varepsilon)|\nabla v|^2 +\Sigma(t,x,v)] . 
\end{align}
\end{lemma}

\begin{proof}
We start with \eqref{R pq phi equality-2} in Lemma \ref{evolution Rpq-2.n}. 
Referring to the second and third lines on the right-hand side of the last equation (the expression with coefficient 
$\zeta(t) |\nabla v|^{s-2}/v^q$), suppressing $\zeta$ for now and denoting the remaining expression by ${\bf I}$ we have
\begin{align}\label{eq-10.17-2}
{\bf I}=& \frac{|\nabla v|^{s-2}}{v^q} \bigg\{ s (p-1)\left [|\nabla v|^2 \Delta_f v - v |\nabla \nabla v|^2
- v \frac{\langle \nabla f , \nabla v \rangle^2}{(m-n)}\right]\nonumber\\
&+ 2s[q(p-1)+1-\varepsilon/2] \nabla \nabla v (\nabla v, \nabla v) - q[q(p-1) +p-\varepsilon] \frac{|\nabla v|^{4}}{v}\bigg\}\nonumber\\
%=&~\frac{|\nabla v|^{s-2}}{v^q} \bigg\{ s (p-1)\left [|\nabla v|^2 \Delta v - v |\nabla \nabla v|^2\right]
%+ 2s[q(p-1)+1-\varepsilon/2] \nabla \nabla v (\nabla v, \nabla v) \bigg\}\nonumber\\
%& -s (p-1) \frac{|\nabla v|^{s-2}}{v^{q-1}} \left[ \frac{|\nabla v|^2}{v} \langle \nabla f, \nabla v\rangle 
%+ \frac{\langle \nabla f , \nabla v \rangle^2}{(m-n)} \right] 
%- q[q(p-1) +p-\varepsilon] \frac{|\nabla v|^{s+2}}{v^{q+1}}\nonumber\\
=& \frac{|\nabla v|^{s-2}}{v^q} \bigg\{ \overbrace{s (p-1)\left [|\nabla v|^2 \Delta v - v |\nabla \nabla v|^2\right]
+ 2s[q(p-1)+1-\varepsilon/2] \nabla \nabla v (\nabla v, \nabla v)}^{{\bf II} 
= s (p-1) [|\nabla v|^2 {\rm tr} A - v |A|^2] + 2s [q(p-1)+1 - \varepsilon/2] |\nabla v|^2 A(e,e) 
} \bigg\} \nonumber\\
& -s (p-1) \frac{|\nabla v|^{s-2}}{v^{q-1}} \left[ \sqrt{m-n}\frac{|\nabla v|^2}{2v} 
+ \frac{\langle \nabla f , \nabla v \rangle}{\sqrt{m-n}} \right]^2\nonumber\\
&+ \left[\frac{s}{4} (p-1) (m-n)-  q[q(p-1) +p-\varepsilon]\right] \frac{|\nabla v|^{s+2}}{v^{q+1}}.
\end{align}

Now for the expression inside the curly brackets on the right-hand side of the second equality, say ${\bf II}$, 
setting $A=\nabla \nabla v$ with ${\rm tr} A=\Delta v$ and $e = \nabla v/ |\nabla v|$, we can write 
\begin{align}
{\bf II} &= s (p-1) \left [|\nabla v|^2 \Delta v - v |\nabla \nabla v|^2\right]
- s\{\varepsilon -2[q(p-1)+1]\} \nabla \nabla v (\nabla v, \nabla v) \\
%=&~ s (p-1) |\nabla v|^2 \frac{{\rm tr} A}{|A|} |A| - s (p-1) v |A|^2 - s\{\varepsilon -2[q(p-1)+1]\} |\nabla v|^2\frac{A(e,e)}{|A|} |A|\nonumber\\
&=-s (p-1) v |A|^2 + \left[ s (p-1)\frac{{\rm tr} A}{|A|} - s\{\varepsilon -2[q(p-1)+1]\} \frac{A(e,e)}{|A|}\right] |\nabla v|^2 |A| \nonumber\\
&=-s (p-1) v \left[ |A|^2 + \frac{|A||\nabla v|^2}{(p-1)v}\left(\{\varepsilon -2[q(p-1)+1]\} \frac{A(e,e)}{|A|}-(p-1)\frac{{\rm tr} A}{|A|} \right) \right], \nonumber
\end{align}
or by extracting a perfect square and rearranging terms 
\begin{align} \label{eq-10.19-2}
{\bf II}=&-s (p-1) v \left[|A| + \frac{1}{(p-1)}\left(\{\varepsilon -2[q(p-1)+1]\} \frac{A(e,e)}{|A|}-(p-1)\frac{{\rm tr} A}{|A|}\right)\frac{|\nabla v|^2}{2v}\right]^2\nonumber\\
& + \frac{sv}{(p-1)} \left[ \{\varepsilon -2[q(p-1)+1]\} \frac{A(e,e)}{|A|}-(p-1)\frac{{\rm tr} A}{|A|}\right]^2 \frac{|\nabla v|^4}{4v^2} \nonumber \\
\overbrace{\le}^{s(p-1)>0}&~\frac{s|\nabla v|^4}{4(p-1)v} \left[ \{\varepsilon -2[q(p-1)+1]\} \frac{A(e,e)}{|A|}-(p-1)\frac{{\rm tr} A}{|A|}\right]^2. 
\end{align}
%where in deducing the inequality on the last line we have made note of $(p-1)v>0$. 
Applying Lemma \ref{lem9.1} to the last expression on the right with the choices of constants $a= \varepsilon -2[q(p-1)+1]$, $b= -(p-1)$ and  
\begin{align}
(a+b)^2 +(n-1) b^2 = \{\varepsilon- 2[q(p-1)+1] -(p-1)\}^2+ (n-1)(p-1)^2,
\end{align}
then gives  
\begin{align}
{\bf II} 
%&= s (p-1) \left [|\nabla v|^2 \Delta v - v |\nabla \nabla v|^2\right]
%- s\{\varepsilon -2[q(p-1)+1]\} \nabla \nabla v (\nabla v, \nabla v)\nonumber\\
&\le \frac{s|\nabla v|^4}{4(p-1)v} \left(\{\varepsilon -2[q(p-1)+1]\} \frac{A(e,e)}{|A|}-(p-1)\frac{{\rm tr} A}{|A|}\right)^2 \nonumber \\
&\le \frac{s|\nabla v|^4}{4(p-1)v} \left[ \{\varepsilon- 2[q(p-1)+1] -(p-1)\}^2+ (n-1)(p-1)^2\right].
\end{align}
As a result using this inequality for ${\bf II}$ and substituting back into ${\bf I}$ in \eqref{eq-10.17-2} gives 
\begin{align}
{\bf I} %=&~\frac{|\nabla v|^{s-2}}{v^q} \bigg\{ s (p-1)\left [|\nabla v|^2 \Delta_f v - v |\nabla \nabla v|^2
%- v \frac{\langle \nabla f , \nabla v \rangle^2}{(m-n)}\right]\\
%&+ 2s[q(p-1)+1-\varepsilon/2] \nabla \nabla v (\nabla v, \nabla v) - q[q(p-1) +p-\varepsilon] \frac{|\nabla v|^{4}}{v}\bigg\}\nonumber\\
\le& \bigg\{ s (p-1)\left [|\nabla v|^2 \Delta v - v |\nabla \nabla v|^2\right]
+ 2s[q(p-1)+1-\varepsilon/2] \nabla \nabla v (\nabla v, \nabla v) \bigg\} \frac{|\nabla v|^{s-2}}{v^q} \nonumber\\
&+ \left[\frac{s}{4} (p-1) (m-n)-  q[q(p-1) +p-\varepsilon]\right] \frac{|\nabla v|^{s+2}}{v^{q+1}} \nonumber\\
\le& \frac{s}{4(p-1)} \left[ \{\varepsilon- 2[q(p-1)+1] -(p-1)\}^2+ (n-1)(p-1)^2\right] \frac{|\nabla v|^{s+2}}{v^{q+1}} \nonumber\\
&+ \left[\frac{s}{4} (p-1) (m-n) - q[q(p-1) +p-\varepsilon]\right] \frac{|\nabla v|^{s+2}}{v^{q+1}} \nonumber\\
\le& \bigg[ \frac{s \{\varepsilon- 2[q(p-1)+1] -(p-1)\}^2+ s(p-1)^2 (m-1)}{4(p-1)} \bigg] \frac{|\nabla v|^{s+2}}{v^{q+1}} \nonumber \\
&- \bigg[ q[q(p-1) +p-\varepsilon]\bigg] \frac{|\nabla v|^{s+2}}{v^{q+1}}. \nonumber 
\end{align} 
Hence from \eqref {R pq phi equality-2} we obtain 
\begin{align}\label{eq-15.25-2}
{\mathscr L}_v^{p}({\mathsf H}^{s,q}_\zeta[v])  &- \varepsilon \zeta(t)\langle \nabla v, \nabla (|\nabla v|^s/v^q) \rangle \nonumber\\
\le &~ \zeta'(t) \frac{|\nabla v|^s}{v^q} - s \zeta (t) \frac{|\nabla v|^{s-2}}{v^q} \left [\frac{1}{2} \partial_t g 
+ (p-1) v {\mathscr Ric}_f^m(g) \right] ( \nabla v, \nabla v) \nonumber\\
&+ \frac{s\zeta (t)}{4(p-1)} \bigg[\{\varepsilon- 2[q(p-1)+1] -(p-1)\}^2 + (p-1)^2(m-1) \nonumber\\
&- \frac{4(p-1)}{s} q[q(p-1) +p-\varepsilon]\bigg] \frac{|\nabla v|^{s+2}}{v^{q+1}}
- \frac{s\zeta(t)}{2v^{q-1}}(p-1) \langle \nabla |\nabla v|^{s-2} , \nabla |\nabla v|^2 \rangle\nonumber\\
&+ \zeta (t) s \frac{|\nabla v|^{s-2}}{v^q}\langle \nabla v, \nabla \Sigma(t,x,v) \rangle
- \zeta (t)\frac{q}{v^{q+1}} |\nabla v|^s\Sigma(t,x,v) \nonumber \\
&+ \Gamma'(v) [|\nabla v|^2 +\Sigma(t,x,v)]
- (p-1)v\Gamma''(v)  |\nabla v|^2. 
\end{align}

Finally by recalling the definition of ${\mathsf H}^{s,q}_\zeta[v]$ in \eqref{Xpq def-2.n} 
we can rewrite \eqref{eq-15.25-2} after a suitable rearrangement as 
\begin{align}
{\mathscr L}_v^{p} ({\mathsf H}^{s,q}_\zeta[v]) 
&- \varepsilon \langle \nabla v, \nabla {\mathsf H}^{s,q}_\zeta[v] \rangle \nonumber\\
\le &~ \zeta'(t) \frac{|\nabla v|^s}{v^q} - s \zeta (t) \frac{|\nabla v|^{s-2}}{v^q} 
\left[ \frac{1}{2} \partial_t g + (p-1) v {\mathscr Ric}_f^m(g) \right] ( \nabla v, \nabla v) \nonumber\\
&+ \frac{s\zeta (t)}{4(p-1)} \bigg[\{\varepsilon- 2[q(p-1)+1] -(p-1)\}^2 + (p-1)^2(m-1) \nonumber\\
&- \frac{4(p-1)}{s} q[q(p-1) +p-\varepsilon]\bigg] \frac{|\nabla v|^{s+2}}{v^{q+1}}
- \frac{s\zeta(t)}{2v^{q-1}}(p-1) \langle \nabla |\nabla v|^{s-2} , \nabla |\nabla v|^2 \rangle\nonumber\\
&+ \zeta (t) s \frac{|\nabla v|^{s-2}}{v^q} \langle \nabla v, \nabla \Sigma(t,x,v) \rangle
- \zeta (t)\frac{q}{v^{q+1}} |\nabla v|^s\Sigma(t,x,v) \nonumber \\
&+ \Gamma'(v) [(1-\varepsilon) |\nabla v|^2 +\Sigma(t,x,v)]
- (p-1)v\Gamma''(v)  |\nabla v|^2. 
\end{align}

An application of super flow \eqref{eq11c-2.n} along with 
$\langle \nabla v, \nabla \Sigma \rangle = \langle \nabla v, \Sigma_x \rangle + \Sigma_v |\nabla v|^2$ 
and the inequality $\langle \nabla |\nabla v|^{s-2}, \nabla |\nabla v|^2\rangle \ge 0$ 
gives the desired conclusion. Note that the last inequality follows by writing $V=|\nabla v|^2$, setting 
$\alpha=(s-2)/2$ (recall $s\ge2$) and observing that 
$\langle \nabla V^\alpha, \nabla V \rangle = \alpha V^{\alpha-1} |\nabla V|^2 \ge 0$.
\end{proof}

 \begin{lemma} \label{evolution Rpq-2.n-inequality-2}
Under the assumptions of Lemma $\ref{evolution Rpq-2.n}$, if the metric $g$ and potential $f$ evolve under 
the flow inequality \eqref{eq11c-2.n} then ${\mathsf H} = \zeta(t) v^{s/[2(s-1)(p-1)]} |\nabla v|^s + \Gamma(v)$ 
satisfies the evolution inequality 
\begin{align} \label{EQ-eq-4.4-2.n-2}
{\mathscr L}_v^{p} ({\mathsf H}) - & p \langle \nabla v, \nabla {\mathsf H} \rangle \nonumber\\
\le & \left\{ \zeta'(t) + s \mathsf{k} \zeta(t) + s \zeta(t) \left[ \Sigma_v(t,x,v)]
+ \frac{\Sigma(t,x,v)/v}{2(s-1)(p-1)} \right] \right\} v^\frac{s}{2(s-1)(p-1)} |\nabla v|^{s} \nonumber\\
&+ s\zeta (t) \bigg[\frac{(s-1)(m-1)(p-1)^2-1}{4(s-1)(p-1)}\bigg] v^{-1+\frac{s}{2(s-1)(p-1)}} |\nabla v|^{s+2} \nonumber\\
&+ s \zeta (t) \langle \nabla v, \Sigma_x(t,x,v) \rangle v^\frac{s}{2(s-1)(p-1)} |\nabla v|^{s-2}
+ \Gamma'(v) \Sigma(t,x,v) \nonumber\\
&- (p-1) [\Gamma'(v)+v\Gamma''(v)] |\nabla v|^2. 
\end{align}
Here ${\mathsf H} = {\mathsf H}^{s,q}_\zeta[v]$ is as in \eqref{Xpq def-2.n} with $q=q^\star = - s/[2(s-1)(p-1)]$ while $s \ge 2$.
\end{lemma}

\begin{proof}
We start with inequality \eqref{EQ-eq-4.4-2.n} established in Lemma \ref{evolution Rpq-2.n-inequality}. Setting $\Omega (q, \varepsilon)$ 
to be the quadratic function 
\begin{align}
\Omega (q, \varepsilon) =&~ \{\varepsilon- 2[q(p-1)+1] -(p-1)\}^2 + (p-1)^2(m-1) \nonumber\\
&- \frac{4(p-1)}{s} q[q(p-1) +p-\varepsilon],
\end{align}
the latter inequality can be written as 
\begin{align} \label{EQ-eq-4.4-2.n-inside}
{\mathscr L}_v^{p} ({\mathsf H}^{s,q}_\zeta[v]) 
- &~\varepsilon \langle \nabla v, \nabla {\mathsf H}^{s,q}_\zeta[v] \rangle \nonumber\\
\le& \left\{ \zeta'(t) + s \mathsf{k} \zeta(t) + \zeta(t) \left[ s \Sigma_v(t,x,v)
- q \frac{\Sigma(t,x,v)}{v} \right] \right\} \frac{|\nabla v|^s}{v^q} \nonumber\\
&+ \frac{s\Omega(q,\varepsilon)\zeta (t)}{4(p-1)} \frac{|\nabla v|^{s+2}}{v^{q+1}}
+ s \zeta (t) \frac{|\nabla v|^{s-2}}{v^q}\langle \nabla v, \Sigma_x (t,x,v)\rangle  \nonumber\\
&- (p-1)v\Gamma''(v)  |\nabla v|^2 + \Gamma'(v) [(1-\varepsilon)|\nabla v|^2 +\Sigma(t,x,v)]. 
\end{align}
A straightforward optimisation of $\Omega$ leads to the optimiser $(q^\star, \varepsilon^\star)$ with 
\begin{align}\label{EQ-4.28}
q^\star = - \frac{s}{2(s-1)(p-1)} , \qquad \text{and} \qquad \varepsilon^\star = p,
\end{align}
and subsequently the optimal value 
\begin{align}\label{EQ-4.29}
\Omega (q^\star, \varepsilon^\star) =-\frac{1}{s-1} +(p-1)^2(m-1).
\end{align}
Substituting the values \eqref{EQ-4.28} and \eqref{EQ-4.29} back into \eqref{EQ-eq-4.4-2.n-inside} gives the desired inequality
\end{proof}

\begin{remark}
The value $\Omega(q^\star, \varepsilon^\star)$ is non-positive {\it iff} $1<p<1+1/\sqrt{(s-1)(m-1)}$ as can be easily seen 
from \eqref{EQ-4.29}. The importance of this non-positive sign will become apparent in the next lemmas where we apply 
maximum principle. 
\end{remark}

\begin{corollary} \label{first-cor-last}
Let $(\mathscr M,g,e^{-f} dv_g)$ be a smooth metric measure space with $\mathscr M$ closed. 
Let $u$ be a positive smooth solution to \eqref{eq11} where $1<p \le 1+1/\sqrt{(s-1)(m-1)}$ with $s \ge 2$ and set $v =pu^{p-1}/(p-1)$. 
Assume the metric $g$ and potential $f$ satisfy the super flow inequality \eqref {eq11c-2.n}. Moreover suppose 
the following conditions hold for some $a$: 
\begin{itemize}
\item $\Gamma'(v) \Sigma(t,x,v) \le 0$,
\item $\Gamma'(v) + v \Gamma''(v) \ge 0$,
\item $\langle \nabla v, \Sigma_x(t,x,v) \rangle \le 0$,  
\item $\Sigma_v(t,x,v) + \Sigma(t,x,v)/[2(s-1)(p-1)v] \le a$.
\end{itemize}
Then for all $x \in \mathscr M$ and $0<t \le T$ we have 
\begin{equation}
[v^\frac{s}{2(s-1)(p-1)} |\nabla v|^s] (x,t) 
\le e^{s(\mathsf{k}+a)t} \left\{ \max_{\mathscr M} 
\left[ v^\frac{s}{2(s-1)(p-1)} |\nabla v|^s + \Gamma(v) \right]_{t=0} - \Gamma(v(x,t)) \right\}. 
\end{equation}
\end{corollary}

\begin{proof}
Using \eqref{EQ-eq-4.4-2.n-2} and the assumptions on $\Gamma$ and $\Sigma$ in the statement of the corollary, we can write for 
${\mathsf H} = \zeta(t) v^{s/[2(s-1)(p-1)]} |\nabla v|^s + \Gamma(v)$
\begin{align} \label{Eq-7.16}
{\mathscr L}_v^{p} & ({\mathsf H}) 
- p \langle \nabla v, \nabla {\mathsf H} \rangle \nonumber\\
\le & \left\{ \zeta'(t) + s\mathsf{k} \zeta(t) + s \zeta(t) \left[ \Sigma_v(t,x,v)]
+ \frac{\Sigma(t,x,v)/v}{2(s-1)(p-1)}\right] \right\} v^\frac{s}{2(s-1)(p-1)} |\nabla v|^{s} \nonumber\\
&+ s\zeta (t) \bigg[\frac{(s-1)(m-1)(p-1)^2-1}{4(s-1)(p-1)}\bigg] v^{-1+ \frac{s}{2(s-1)(p-1)}} |\nabla v|^{s+2} \nonumber\\
&+ s \zeta (t) \langle \nabla v, \Sigma_x(t,x,v) \rangle v^\frac{s}{2(s-1)(p-1)} |\nabla v|^{s-2}
+ \Gamma'(v)\Sigma(t,x,v) \nonumber\\
&-(p-1)[\Gamma'(v)+v\Gamma''(v)]  |\nabla v|^2 \nonumber \\
\le &~ \left[\zeta'(t) + s({\mathsf k}+a) \zeta (t)\right] v^\frac{s}{2(s-1)(p-1)}|\nabla v|^{s} \nonumber\\
&+ s\zeta (t) \bigg[\frac{(s-1)(m-1)(p-1)^2-1}{4(s-1)(p-1)}\bigg] v^{-1+ \frac{s}{2(s-1)(p-1)}} |\nabla v|^{s+2}.
\end{align}
The function $\zeta(t) = e^{-s(\mathsf{k}+a)t}$ is non-negative, smooth and satisfies $\zeta'+s(\mathsf{k}+a)\zeta=0$. 
Thus substituting in \eqref{Eq-7.16} and noting $(s-1)(m-1)(p-1)^2 \le 1$ due to the range of $p$ we have 
\begin{equation}
{\mathscr L}_v^{p} ({\mathsf H}) - p \langle \nabla v, \nabla {\mathsf H} \rangle \le 0.
\end{equation} 
The assertion is now a consequence of the weak maximum principle giving 
\begin{align}
{\mathsf H}(x,t) &= e^{-s(\mathsf{k}+a)t} [v^\frac {s}{2(s-1)(p-1)} |\nabla v|^s](x,t) + \Gamma(v(x,t)) \nonumber \\
&\le \max_\mathscr M [v^\frac{s}{2(s-1)(p-1)} |\nabla v|^s + \Gamma(v)]_{t=0} = \max_\mathscr M {\mathsf H}(x,0).
\end{align}
The proof is thus complete. 
\end{proof}

\begin{corollary} \label{second-cor-last}
Let $(\mathscr M,g,e^{-f} dv_g)$ be a smooth metric measure space with $\mathscr M$ closed. 
Let $u$ be a positive smooth solution to \eqref{eq11} where $1<p \le 1+1/\sqrt{m-1}$ and set $v =pu^{p-1}/(p-1)$. 
Assume the metric $g$ and potential $f$ satisfy the super flow inequality \eqref {eq11c-2.n} with $\mathsf{k} \ge 0$. Moreover suppose 
the following conditions hold: 
\begin{itemize}
\item $\Sigma(t,x,v) \le 0$,
\item $\langle \nabla v, \Sigma_x(t,x,v) \rangle \le 0$, 
\item $\Sigma_v(t,x,v) + \Sigma(t,x,v)/[2(p-1)v] \le 0$. 
\end{itemize}
Then for $x \in \mathscr M$ and $0<t \le T$ we have 
\begin{equation}
\frac{p^2t}{1+2{\mathsf k}t} [v^\frac{1}{p-1} |\nabla v|^2](x,t) \le (p-1) \left[ \max_\mathscr M v^\frac{p}{p-1}(x,0) - v^\frac{p}{p-1}(x,t) \right]. 
\end{equation}
\end{corollary}

\begin{proof}
For ${\mathsf H}[v] = \zeta(t) v^{1/(p-1)} |\nabla v|^2 + (p-1)v^{p/(p-1)}/p^2$ where 
$\Gamma(x,v)=(p-1)v^{p/(p-1)}/p^2$, $s=2$ and $q^\star=-1/(p-1)$, we have from 
Lemma \ref{evolution Rpq-2.n-inequality-2},
\begin{align} 
{\mathscr L}_v^{p} ({\mathsf H}) 
-&~p \langle \nabla v, \nabla {\mathsf H} \rangle \nonumber\\
\le & \left\{ \zeta'(t) + 2\mathsf{k} \zeta (t) + 2 \zeta(t)  \left[\Sigma_v(t,x,v)]
+ \frac{\Sigma(t,x,v)}{2(p-1)v} \right] \right\} v^\frac{1}{p-1} |\nabla v|^2 \nonumber\\
&+ \frac{(m-1)(p-1)^2-1}{2(p-1)} \zeta(t) v^\frac{2-p}{p-1} |\nabla v|^4 
+ 2 \zeta (t) \langle \nabla v, \Sigma_x(t,x,v) \rangle v^\frac{1}{p-1} \nonumber \\
&+ \Gamma'(v) \Sigma(t,x,v) -(p-1) [\Gamma'(v)+v\Gamma''(v)] |\nabla v|^2 \nonumber \\
\le& \left\{ \zeta'(t) + 2\mathsf{k} \zeta (t) + 2 \zeta(t)  \left[\Sigma_v(t,x,v)]
+ \frac{\Sigma(t,x,v)}{2(p-1)v} \right] \right\} v^\frac{1}{p-1} |\nabla v|^2 \nonumber\\
&+ \frac{(m-1)(p-1)^2-1}{2(p-1)} \zeta(t) v^\frac{2-p}{p-1} |\nabla v|^4 
+ 2 \zeta (t) \langle \nabla v, \Sigma_x(t,x,v) \rangle v^\frac{1}{p-1} \nonumber \\
&+ \frac{1}{p} v^\frac{1}{p-1} \Sigma(t,x,v) - v^\frac{1}{p-1} |\nabla v|^2.
\end{align}
In particular, with the assumptions in the statement of the corollary this gives 
\begin{align} 
{\mathscr L}_v^{p} ({\mathsf H}) 
- p \langle \nabla v, \nabla {\mathsf H} \rangle 
\le&~[\zeta'(t) + 2\mathsf{k} \zeta(t) -1] v^\frac{1}{p-1} |\nabla v|^2 \nonumber \\
&~+ \frac{(m-1)(p-1)^2-1}{2(p-1)} \zeta(t) v^\frac{2-p}{p-1} |\nabla v|^4.  
\end{align}

Now when $\mathsf{k} \ge 0$ by taking $\zeta(t) = t/(1+2\mathsf{k} t)$ we have 
$(\zeta' + 2 \mathsf{k} \zeta -1) \le 0$. Therefore subject to $1<p \le 1+1/\sqrt{m-1}$ 
we have ${\mathscr L}_v^{p}({\mathsf H}) - p \langle \nabla v, \nabla {\mathsf H} \rangle \le 0$. 
The conclusion now follows by an application of the weak maximum principle.  
\end{proof}

\qquad \\
{\bf Acknowledgement.} The authors gratefully acknowledge support from the Engineering and Physical 
Sciences Research Council (EPSRC) through the grant EP/V027115/1. All data is provided in full in the 
results section. Additional data is in the public domain at locations cited in the reference section.

\end{document}